\documentclass[a4paper,11pt,leqno]{amsart}

\usepackage[a4paper,margin=2.5cm,headsep=0.2in,headheight=14pt]{geometry}

\usepackage{amsmath,amscd,amsfonts,amssymb}

\usepackage{tikz}

\usepackage[english]{babel}

\usepackage{calrsfs}

\usepackage{fancyhdr}

\numberwithin{equation}{section}

\newtheorem{theorem}{Theorem}[section]

\newtheorem{prop}[theorem]{Proposition}

\newtheorem{lemma}[theorem]{Lemma}

\theoremstyle{definition}

\newtheorem{defi}[theorem]{Definition}

\newtheorem{rem}[theorem]{Remark}

\newcommand{\RR}{\mathbb{R}}

\newcommand{\ZZ}{\mathbb{Z}}

\newcommand{\TT}{\mathit{TT}}

\newcommand{\TK}{\mathit{TK}}

\newcommand{\KK}{\mathit{KK}}

\title{Collapsing flat manifolds with totally geodesic boundary}

\author{Chady El Mir}

\dedicatory{In memory of Ahmad El Soufi, a distinguished mathematician and dear friend.}

\address{Lebanese University, Laboratoire de Math\'ematiques et Applications (LaMA), Tripoli, Lebanon}

\email{chady.mir@gmail.com}

\graphicspath{{images/}}


\begin{document}

\date{}

\selectlanguage{english}

\begin{abstract}

We study volume collapse of compact connected flat manifolds with nonempty totally geodesic boundary and uniformly bounded diameter. We prove subsequential Gromov--Hausdorff convergence to compact flat orbifolds, allowing reflector loci, and describe the limit using the canonical involutions of the doubles. We distinguish the cases in which the maximal distance to the boundary tends to zero or to a positive limit, and classify all possible limits in dimension three.

\end{abstract}

\medskip

\maketitle

\selectlanguage{english}
\par
\medskip

\section{Introduction and results}\label{intro}

\subsection*{Collapse of Riemannian manifolds}

For a fixed dimension $n$, the class of closed connected Riemannian $n$-manifolds with uniformly bounded diameter and a uniform lower Ricci curvature bound is precompact in the Gromov--Hausdorff topology \cite{gromov}. If a sequence in this class converges to $X$, one says that it \emph{collapses} when the Hausdorff dimension of $X$ is smaller than $n$. Under a uniform lower sectional curvature bound, this is equivalent to $\operatorname{Vol}(M_i,g_i)\to0$.

Collapse with uniformly bounded sectional curvature is described by the $F$-structures of Cheeger--Gromov \cite{cheegergromov,cheegergromov2} (see \cite{pansu85} for a survey) and the nilpotent Killing structures of Cheeger--Fukaya--Gromov \cite{cff}, defined for suitable nearby invariant metrics on collapsed regions. When closed manifolds with uniformly bounded sectional curvature converge to a lower-dimensional closed Riemannian manifold, Fukaya's fibration theorem gives, for sufficiently large $i$, a fibre bundle over that limit with infranil fibres \cite{fukaya}. Here the Riemannian limit is fixed and has positive injectivity radius, as required in the fibration theorem.

Under a uniform lower sectional curvature bound, limits of closed Riemannian $n$-manifolds are Alexandrov spaces with the same curvature lower bound \cite{bgp}. Equivariant Gromov--Hausdorff convergence keeps track of group actions as well as the metric spaces: when pointed covers and their deck transformation groups converge equivariantly, their quotients converge to the quotient by the limiting group (Theorem~\ref{equivariant}, following \cite{fukayayama,fukaya1}).

\subsection*{The flat case}

Flat manifolds provide a setting in which collapse can be studied through crystallographic groups. By Bieberbach's theorems \cite{charlap}, a closed flat manifold is $\mathbb R^n/\Gamma$, where $\Gamma$ is a discrete cocompact torsion-free group of Euclidean isometries and its translation subgroup is a lattice of finite index.

The almost-flat theorem of Gromov and Ruh \cite{gromovalmost,ruh} asserts that a closed Riemannian $n$-manifold is diffeomorphic to an infranilmanifold if $\operatorname{diam}(M,g)^2\sup|\sec_g|$ is sufficiently small, with a threshold depending only on $n$. Bounded curvature and bounded diameter together with volume collapse do not imply this smallness condition: $S^2(1)\times S^1(\varepsilon)$, with circle radius $\varepsilon\to0$, collapses to $S^2$ and is not infranil, since its second homotopy group is nonzero. For flat metrics, the following result is the starting point of this paper.

\medskip

\noindent\textbf{The theorem of Bettiol--Derdzinski--Piccione.} Bettiol, Derdzinski and Piccione \cite{betderpic} study the Teichm\"uller and moduli spaces of flat metrics on closed manifolds (for closed flat $3$-manifolds see \cite{kangkim,kang3d}). Their Theorem~A shows that a Gromov--Hausdorff limit of closed connected flat manifolds of a fixed dimension is a compact flat orbifold, and that every compact connected flat orbifold is obtained as such a limit in a suitable fixed dimension. In particular, a sequence of closed flat $n$-manifolds with uniformly bounded diameter has a subsequence converging to a compact flat orbifold.

In dimensions at most three the possible collapsed limits are explicitly classified. A sequence of closed flat surfaces can limit only to a point, a circle, an interval, a flat torus or a flat Klein bottle; in the strictly collapsing case only the first three occur. A flat $3$-torus collapses only to a point, a circle or a flat $2$-torus. By \cite[Theorem~D]{betderpic}, corrected as in Remark~\ref{rem:bdp}, the two-dimensional collapsed limits of closed flat $3$-manifolds are the ten orbifolds $\circ$, $\times\times$, $**$, $*\times$, $2222$, $333$, $442$, $632$, $22{*}$ and $22{\times}$. These results provide the general orbifold-limit theorem and the low-dimensional classification used below.

\medskip

For manifolds with totally geodesic boundary, collapse must also account for the boundary and the involution of the double. Reflector loci can already occur in limits of closed flat manifolds: a flat Klein bottle can collapse to an interval, and its product with a fixed circle can collapse to an annulus with mirror boundary \cite[Section~6]{betderpic}.

A separate issue is the distinction between intrinsic boundary distances and distances restricted from the ambient manifold. For example, the flat M\"obius band $M_w=(\mathbb R\times[-w,w])/\langle(x,y)\mapsto(x+1,-y)\rangle$ collapses, as $w\to0$, to a circle of length~$1$, while its boundary with the intrinsic metric is a circle of length~$2$.

\subsection*{Results}

In this paper we study sequences of compact connected flat manifolds with nonempty \emph{totally geodesic} boundary and uniformly
bounded diameter. (In the flat case, the double of such a manifold along its boundary is a closed flat manifold, which allows one
to apply the results of \cite{betderpic}.) Let $M^n$ be a compact connected manifold with nonempty boundary,
and let $g_i$ be flat metrics on $M$ with totally geodesic boundary, $\operatorname{diam}(M,g_i)\le D$ and $\operatorname{Vol}(M,g_i)\to0$. Our first result
describes the limit by combining the theorem of \cite{betderpic} for the doubles with the equivariant Gromov--Hausdorff convergence
of \cite{fukayayama}. The quantity that distinguishes the two cases is the largest distance to the boundary,
$\tau_i=\sup_{x\in M}d_{g_i}(x,\partial M)$. In the terminology of Yamaguchi and Zhang \cite{yz}, $\tau_i$ is the inradius of $(M,g_i)$, and the case $\tau_\infty=0$ of Theorem~A below is an inradius collapse. Convergence and collapse of manifolds with boundary under lower curvature bounds have been studied by Yamaguchi and Zhang \cite{yz,yz1,yz2} (see also \cite{wong}), and under lower Ricci bounds by Huang and Yamaguchi \cite{hy}. Flat metrics with totally geodesic boundary satisfy the curvature hypotheses of these works. In particular, for inradius collapse, Huang and Yamaguchi \cite{hy} show that the limit of the boundaries carries an isometric involution whose quotient is the limit of the manifolds, and that the boundary has at most two components. In the flat case, Theorem~A(ii) and Theorem~\ref{theorem:slab}(4) recover these facts with an explicit description of the involution, and identify the limit as a compact flat orbifold; the boundary bound holds for every compact flat manifold with totally geodesic boundary, not only in the collapsed case.

\medskip

\noindent\textbf{Theorem A.}

\emph{Let $M^n$, $n\ge2$, be a compact, connected manifold with nonempty boundary. Suppose there is a sequence of flat Riemannian metrics $g_i$ on $M$ with totally geodesic boundary such that $\operatorname{diam}(M,g_i)\le D$ for some fixed $D>0$ and $\operatorname{Vol}(M,g_i)\to0$. Put $\tau_i=\sup_{x\in M}d_{g_i}(x,\partial M)$. Then, after passing to a subsequence, $\tau_i\to\tau_\infty\in[0,D]$ and $(M,g_i)$ converges in the Gromov--Hausdorff sense to a compact flat orbifold $X$ (flat orbifolds are allowed to have reflector loci). Moreover:}

\begin{itemize}

\item[\textup{(i)}] \emph{(Non-degenerate case).} If $\tau_\infty>0$, then, after passing to a further subsequence, the doubles $\widetilde M_i=M\cup_{\partial M}M$ converge to a compact flat orbifold $O$, the doubling involutions $\iota_i$ converge in the equivariant Gromov--Hausdorff sense to a non-trivial isometric involution $\bar\iota\in\mathrm{Isom}(O)$, and $X=O/\langle\bar\iota\rangle$.

\item[\textup{(ii)}] \emph{(Degenerate case).} If $\tau_\infty=0$, then $\tau_i\to0$ and:

\begin{itemize}

\item[\textup{(a)}] if $\partial M$ is disconnected, then it has exactly two components, and $X=\lim_{GH}(\Sigma_i,g_i|_{\Sigma_i})$ for either boundary component $\Sigma_i$ of $(M,g_i)$, with its induced intrinsic metric;

\item[\textup{(b)}] if $\partial M$ is connected, then $(\partial M,g_i|_{\partial M})$ carries a free isometric involution $\sigma_i$ (the deck involution of the Riemannian double cover $\partial M\to N_i$ given by Theorem~\ref{theorem:slab}(4b) for the metric $g_i$, where $N_i$ is the closed flat $(n-1)$-manifold of Theorem~\ref{theorem:slab}(2)), and 
\[
X=\lim_{GH}\bigl((\partial M,g_i|_{\partial M})/\langle\sigma_i\rangle\bigr),
\]
where each quotient carries its quotient metric.

\end{itemize}

\end{itemize}

\medskip

The hypothesis $\operatorname{Vol}(M,g_i)\to0$ is not used in the proof: it only ensures that the sequence genuinely collapses, and the conclusions of Theorem~A, with the same proof, hold for any sequence of flat metrics with totally geodesic boundary and $\operatorname{diam}(M,g_i)\le D$, in which case the limit $X$ may have dimension~$n$. In the degenerate case, $X$ is also the limit of the doubles of $M_i$. In case (b), the limit of $(\partial M,g_i|_{\partial M})$ itself may differ from $X$: for the M\"obius bands above, $X$ is a circle of length~$1$ while $\partial M_i$ is a circle of length~$2$. A key ingredient is a structure theorem for compact flat manifolds with totally geodesic boundary (Theorem~\ref{theorem:slab}): their universal cover is a slab $\mathbb R^{n-1}\times[-\tau,\tau]$, and in particular the boundary has at most two components.

Our second result gives an explicit list of the possible limits in dimension three. By Section~\ref{klein}, there are, up to isometry
and up to the choice of parameters, five families of compact flat $3$-manifolds with totally geodesic boundary,
namely $T^2\times[-d,d]$, $K\times[-d,d]$ (with disconnected boundary) and the three ``solid M\"obius bands'' $\TT$, $\TK$ and $\KK$
(with connected boundary) \cite{elaf}.

\medskip

\noindent\textbf{Theorem B.}

{\itshape Let $(M_i,g_i)$ be a sequence of compact connected flat $3$-manifolds with nonempty totally geodesic boundary, with
$\operatorname{diam}(M_i,g_i)\le D$ and $\operatorname{Vol}(M_i,g_i)\to0$, and put
$\tau_i=\sup_{x\in M_i}d_{g_i}(x,\partial M_i)$. After passing to a subsequence with $\tau_i\to\tau_\infty$,
$(M_i,g_i)$ converges in the Gromov--Hausdorff sense to a compact flat orbifold $X$ of dimension
$k\in\{0,1,2\}$ (orbifolds are allowed to have reflector loci), and:}

\begin{itemize}

\item[\textup{(i)}] {\itshape if $\tau_\infty>0$, then $X=O/\langle\bar\iota\rangle$, where $O$ is a limit of the
doubles of $M_i$ and $\bar\iota$ is a non-trivial reflection of $O$ (that is, writing $O=\mathbb R^k/\Gamma$, $\bar\iota$ is induced by a Euclidean reflection of $\mathbb R^k$ normalizing $\Gamma$), and $X$ is one of

\begin{itemize}

\item a closed interval $I$;

\item a flat annulus $S^1\times I$ whose two boundary circles are mirrors (the orbifold $**$);

\item a flat rectangle $I\times I$ whose four sides are mirrors, with four corner reflectors of order $2$
(the orbifold $*2222$);

\item a flat M\"obius band whose boundary circle is a mirror (the orbifold $*\times$);

\item the orbifold $2{*}22$: a disc with one cone point of order $2$ and mirror boundary with two corner
reflectors of order $2$;

\item the orbifold $22{*}$: a disc with two cone points of order $2$ and mirror boundary, without corner reflectors;

\end{itemize}}

\item[\textup{(ii)}] {\itshape if $\tau_\infty=0$, then $X$ is one of

\begin{itemize}

\item a point;

\item a circle $S^1$ or a closed interval $I$;

\item a flat torus $T^2$ or a flat Klein bottle $K$.

\end{itemize}}

\end{itemize}

{\itshape Moreover, each of these limits occurs. More precisely, the limits that arise from each of the five families
are listed in Section~\textup{\ref{theoB}}.}

\medskip

Theorem~B exhibits a phenomenon that first appears in dimension three: the rectangle $*2222$ and the orbifold $2{*}22$ occur as limits of compact flat $3$-manifolds with totally geodesic boundary, but neither occurs as a collapsed limit of closed flat $3$-manifolds \cite[Theorem~D]{betderpic}. The orbifold $22{*}$ does occur in the closed case. When $\tau_\infty>0$, the boundary images converge to the fixed set of $\bar\iota$ in the limiting double, as made precise in Remark~\ref{rem:boundary_limits}, and their image in $X$ lies in the reflector locus. When $\tau_\infty=0$, the limit can retain the full two-dimensional soul. Two-dimensional limits with vanishing three-volume also occur in the closed case, for example $T^2\times S^1(\varepsilon)\to T^2$ as $\varepsilon\to0$.

\begin{rem}\label{rem:bdp}
Some computations in \cite[Section~6]{betderpic} need correction. In Examples~6.2 and~6.3 the holonomy acts on $W_2$ by rotations, so the limits $T^2/\mathbb Z_4$ and $T^2/\mathbb Z_6$ are the orientable orbifolds $442$ and $632$, not $4{*}2$ and $3{*}3$. In fact $4{*}2$ and $3{*}3$ cannot occur: by \cite[Proposition~3.2]{betderpic} the point group of a limit is a quotient of the holonomy group, which is abelian for closed flat $3$-manifolds, while the point groups $D_4$ and $D_3$ of $4{*}2$ and $3{*}3$ are not. In Example~6.10 ($G_2$), collapsing $W_1$ gives $\mathbb R^2/(\Lambda\rtimes\{\pm\mathrm{Id}\})$, with $\Lambda=\mathbb Z v_2\oplus\mathbb Z v_3$, which is the pillowcase $2222$, not $22{*}$; this agrees with the case $a_\infty=0<b_\infty$ of $\TK$ in Section~\ref{theoB}. In Example~6.7(a) ($B_4$), the image of $(E,\frac12(v_2+v_3))$ is the glide reflection $(x,y)\mapsto(x+\frac a2,-y+\frac b2)$, and the limit is $22{\times}$, not $2222$. Since $2222$ arises from $G_2$, $22{\times}$ from $G_6$ and $22{*}$ from $B_3$, the remaining entries of \cite[Theorem~D]{betderpic} are unchanged. Finally, $*2222$ and $2{*}22$ do not occur: their point group $\mathbb Z_2^2$ forces the holonomy to be $\mathbb Z_2^2$, so the manifold is $G_6$, $B_3$ or $B_4$, whose lattices are rectangular with axes the three isotypic lines, and collapsing one axis gives $22{\times}$ (all three axes) for $G_6$, $22{*}$, $**$, $\times\times$ for $B_3$, and $22{\times}$, $*\times$, $\times\times$ for $B_4$.
\end{rem}

\medskip

The proofs combine geometric and group-theoretic arguments. For Theorem~A we pass to the double
$\widetilde M_i=M_i\cup_{\partial M_i}M_i$, a closed flat manifold with an isometric involution $\iota_i$ with $\widetilde M_i/\langle\iota_i\rangle=M_i$.

We apply \cite{betderpic} and the equivariant compactness theorem of Fukaya--Yamaguchi to the pairs $(\widetilde M_i,\langle\iota_i\rangle)$,
and we show that the limit involution is trivial exactly when $\tau_\infty=0$. In the degenerate case the slab structure
theorem identifies $X$ with the limit of the soul, a flat manifold of dimension $n-1$. For Theorem~B we write the
doubles as explicit quotients $\mathbb R^3/\widetilde\Gamma_i$ and compute the limits of the groups $\widetilde\Gamma_i$ along each of the five families.

\medskip

The paper is organized as follows. In Section~\ref{flat} we recall facts about flat manifolds and orbifolds, and prove the slab
structure theorem for flat manifolds with totally geodesic boundary. In Section~\ref{theoremA} we recall the equivariant Gromov--Hausdorff convergence
and prove Theorem~A. In Section~\ref{klein} we describe the compact flat $3$-manifolds with totally geodesic boundary, following \cite{elaf}. In Section~\ref{theoB} we study their doubles and prove Theorem~B, by analysing the collapse of each of the five families.

\section{Flat manifolds with geodesic boundary and flat orbifolds}\label{flat}

Closed flat manifolds are quotients $\mathbb{R}^n/\Gamma$, where $\Gamma$ is a discrete cocompact subgroup of affine isometries of $\RR^n$ acting freely. By the theorem of Bieberbach (\cite{bieb1,bieb2}) $\Gamma$  is an extension of a finite group $G$ by a lattice $\Lambda$ of $\RR^n$ (see also \cite{charlap}). This lattice is the subgroup of the elements of $\Gamma$ that are translations, we obtain then the following exact sequence:

$$
0\longrightarrow\Lambda\longrightarrow\Gamma\longrightarrow G\longrightarrow 1
$$

Actually, if $M$ is a closed flat manifold, $M$ is the quotient of the flat torus $\mathbb{R}^n/\Lambda$ by an isometry group isomorphic to $G$.

Two closed flat manifolds $\RR^n/\Gamma$ and $\RR^n/\Gamma^\prime$ are homeomorphic if and only if $\Gamma$ and $\Gamma^\prime$ are isomorphic. In that case the groups are conjugate by an affine isomorphism of $\RR^n$, so the manifolds are affinely diffeomorphic. For open problems on crystallographic groups and Bieberbach groups see \cite{szczepanski1,szczepanski2}.

 The classification of closed flat manifolds of dimension $3$ results from a direct method of classification of discrete, cocompact subgroups of $\mathrm{Isom}(\RR^3)$ operating freely. 

This classification is due to W. Hantzsche and H. Wendt (1935) (see \cite{hantzschewendt}), and exposed in the book \cite{wolfspaces} of J.A. Wolf.

There are ten closed flat $3$-manifolds up to diffeomorphism: six orientable and four non-orientable. In Section~\ref{theoB} we determine which types occur as doubles of compact flat $3$-manifolds with totally geodesic boundary; see also \cite[Section~3.3]{elaf}.
\medskip

\subsection{Flat manifolds with geodesic boundary}

Let $(M,g)$ be a compact connected flat manifold with non-empty totally geodesic boundary, and let $\widetilde M=M\cup_{\partial M}M$ be its double. Near $\partial M$ the metric has the product form $dt^2+g_{\partial M}$ in Fermi coordinates (the shape operators $S_t$ of the equidistant hypersurfaces satisfy $S_t'+S_t^2=0$ by flatness, and $S_0=0$ by total geodesicity, so $S_t\equiv0$). Hence the glued metric $\widetilde g$ is smooth and flat, $\widetilde M$ is a closed connected flat
manifold, and the reflection $\iota$ is a smooth isometric involution with $\mathrm{Fix}(\iota)=\partial M$.

\begin{theorem}[Slab structure]\label{theorem:slab}

Let $(M^n,g)$, $n\ge2$, be a compact connected flat manifold with non-empty totally geodesic
boundary, and let $\tau=\sup_{x\in M}d(x,\partial M)$. Then $\tau>0$ and:

\begin{enumerate}

\item[\textup{(1)}] The universal cover of $M$ is isometric to the slab
$S=H\times[-\tau,\tau]\subset\mathbb R^n$, $H=\mathbb R^{n-1}$, with the Euclidean metric.

The deck group $\Gamma\cong\pi_1(M)$ acts by $\gamma(h,t)=(A_\gamma h,\,s_\gamma t)$, where
$A_\gamma\in\mathrm{Isom}(H)$ and $s_\gamma\in\{\pm1\}$; the maps $\gamma\mapsto A_\gamma$ and
$s:\Gamma\to\{\pm1\}$ are homomorphisms.

\item[\textup{(2)}] The map $\gamma\mapsto A_\gamma$ is injective, and $\bar\Gamma=\{A_\gamma\}$
acts freely and properly discontinuously on $H$. Thus $N:=H/\bar\Gamma$ is a closed flat
$(n-1)$-manifold.

\item[\textup{(3)}] The projection $S\to H$ induces a surjective map $p:M\to N$ with

\[
d_N(p(x),p(y))\le d_M(x,y)\le d_N(p(x),p(y))+2\tau\qquad(x,y\in M).
\]

\item[\textup{(4)}] $\partial M$ has at most two connected components.

\begin{itemize}

\item[\textup{(a)}] If $s\equiv1$, then $\partial M$ has exactly two components, and each of them,
with the intrinsic metric induced by $g$, is isometric to $N$.

\item[\textup{(b)}] If $s\not\equiv1$, then $\partial M$ is connected and
$\partial M\cong H/\bar\Gamma^0$, where $\Gamma^0=\ker s$ has index $2$ in $\Gamma$ and
$\bar\Gamma^0=\{A_\gamma:\gamma\in\Gamma^0\}$. The covering $\partial M\to N$ is a Riemannian
double cover, whose non-trivial deck transformation $\sigma$ is a free isometric involution of
$(\partial M,g|_{\partial M})$ with $(\partial M)/\sigma\cong N$. Moreover, for $x\in\partial M$,
$\sigma(x)$ is the point where the geodesic starting at $x$ orthogonally to $\partial M$
(into $M$) first returns to $\partial M$; this geodesic has length $2\tau$.

\end{itemize}

\end{enumerate}

\end{theorem}

\begin{proof}

Let $\widetilde M$ be the double of $M$, with $M$ the first copy. We know that $\widetilde M$ is a closed connected flat manifold and near $\partial M$ its metric is a
product $dt^2+g_{\partial M}$. Hence the universal cover of $\widetilde M$ is $\mathbb R^n$, with
covering map $\pi:\mathbb R^n\to\widetilde M$, a local isometry.

\medskip

\noindent\emph{Step 1: the universal cover of $M$ is a convex subset $C\subset\mathbb R^n$.}

The set $\pi^{-1}(M)$ is closed, and $\pi:\pi^{-1}(M)\to M$ is a covering map. Since
$\pi^{-1}(M)$ is a manifold with boundary, it is locally connected, so its connected components
are closed. Let $C$ be one of them. Since $M$ is connected and locally path connected,
$\pi|_C:C\to M$ is a covering map; in particular it is surjective.

Let $x\in C$. Choose a small Euclidean ball $U$ around $x$ on which $\pi$ is an isometry onto its
image and, if $\pi(x)\in\partial M$, such that $\pi^{-1}(\partial M)\cap U$ is a piece of an affine
hyperplane (this uses the product structure near $\partial M$, and total geodesicity). Then
$\pi^{-1}(M)\cap U$ is either all of $U$ (if $\pi(x)\notin\partial M$) or the part of $U$ on one
side of that hyperplane. In both cases it is convex, hence connected, so it is contained in $C$.

Thus $C\cap U$ is convex: $C$ is closed, connected and locally convex. By the
Tietze--Nakajima theorem \cite{tietze,nakajima} (a closed, connected, locally convex subset of $\mathbb R^n$ is convex),
$C$ is convex. In particular $C$ is contractible, so $\pi|_C$ is the universal covering of $M$,
and $\Gamma=\pi_1(M)$ is its deck group, acting on $C$ by isometries. Since $C$ is convex, its
induced length metric is the restriction of the Euclidean metric, so $\pi|_C$ is a Riemannian
covering and the deck transformations are restrictions of Euclidean isometries. Also $C$ has
non-empty interior, and $\partial C=C\cap\pi^{-1}(\partial M)\neq\emptyset$, because $\pi|_C$ is
surjective and $\partial M\neq\emptyset$.

\medskip

\noindent\emph{Step 2: $C$ is a slab.}

By Step 1, $\partial C$ is a closed embedded hypersurface which is locally an open piece of an
affine hyperplane. Let $Q$ be a connected component of $\partial C$. The set of points of $Q$ near
which $Q$ coincides with an open piece of a fixed hyperplane $P$ is open in $Q$, and it is closed
in $Q$ (two hyperplanes sharing an open piece coincide), so $Q\subset P$ for a single hyperplane
$P$; moreover $Q$ is open in $P$ and closed in $\mathbb R^n$, so $Q=P$. Hence every component of
$\partial C$ is an affine hyperplane, and distinct components are disjoint, hence parallel.

Let $P\subset\partial C$ be such a component. If $L$ is a supporting hyperplane of $C$ at a point
of $P$, then $P\subset C$ lies on one side of $L$ and touches $L$, so $P=L$. Hence $C$ lies in
one closed half-space bounded by $P$.

There are at most two components: if $P_1,P_2,P_3\subset\partial C$ were three parallel
hyperplanes with $P_2$ between $P_1$ and $P_3$, convexity would give that $C$ contains the slab
between $P_1$ and $P_3$, so $P_2$ would lie in the interior of $C$, contradicting
$P_2\subset\partial C$.

If $\partial C=P$ is a single hyperplane, let $H^+$ be the closed half-space containing $C$. The set
$C\cap\mathrm{int}\,H^+$ is closed in $\mathrm{int}\,H^+$; it is open, since it contains no point of
$\partial C$; it is non-empty, since $C$ has interior points; and $\mathrm{int}\,H^+$ is connected.

So $C\supset\mathrm{int}\,H^+$, and as $C$ is closed, $C=H^+$. Likewise, if $\partial C=P_1\sqcup P_2$,
then $C$ lies between $P_1$ and $P_2$ (each $P_i$ lies in $C$, so $C$ is on the side of $P_i$
containing the other plane), and the same argument shows that $C$ is the closed slab bounded by
$P_1$ and $P_2$.

A half-space is impossible. For $x\in C$ we have
$\mathrm{dist}_{\mathrm{Eucl}}(x,\partial C)=d_M(\pi(x),\partial M)$. Indeed, a minimizing path in
$M$ from $\pi(x)$ to $\partial M$ lifts to a path in $C$ starting at $x$ and ending in
$\pi^{-1}(\partial M)\cap C=\partial C$, which gives $\le$; and the image under $\pi$ of a path in
$C$ from $x$ to $\partial C$ is a path in $M$ from $\pi(x)$ to $\partial M$ of the same length, which
gives $\ge$. Hence $\mathrm{dist}(\cdot,\partial C)\le\tau<\infty$ on $C$, whereas this function is
unbounded on a half-space. So $C$ is a slab. After an isometry of $\mathbb R^n$,
$C=H\times[-w,w]$ with $w>0$, and $\sup_C\mathrm{dist}(\cdot,\partial C)=w$ (attained on
$H\times\{0\}$). Since $\pi|_C$ is onto, $w=\tau$; in particular $\tau>0$. Set $S=C$.

\medskip

\noindent\emph{Step 3: the form of the deck transformations (proof of (1)).}

Let $\gamma\in\Gamma$. Since $S$ has non-empty interior, the isometry $\gamma$ of $S$ extends
uniquely to an affine isometry of $\mathbb R^n$, again denoted $\gamma$. It preserves
$\partial S=H\times\{\pm\tau\}$ and hence permutes these two parallel hyperplanes. Therefore its
linear part preserves the direction space $H\times\{0\}$ and its orthogonal complement
$\{0\}\times\mathbb R$, so $\gamma(h,t)=(A_\gamma h,\theta_\gamma(t))$ with
$A_\gamma\in\mathrm{Isom}(H)$ and $\theta_\gamma$ an isometry of $\mathbb R$ preserving
$\{-\tau,\tau\}$. The only such isometries are $t\mapsto t$ and $t\mapsto-t$, so
$\theta_\gamma(t)=s_\gamma t$ with $s_\gamma=\pm1$ ($s_\gamma=1$ if $\gamma$ preserves each of the
two hyperplanes, $s_\gamma=-1$ if it swaps them). Clearly $\gamma\mapsto A_\gamma$ and $s$ are
homomorphisms.

\medskip

\noindent\emph{Step 4: proof of (2).}

\emph{Freeness.} Let $\gamma\ne\mathrm{id}$ and suppose $A_\gamma h=h$. If $s_\gamma=1$, then
$\gamma$ fixes $(h,t)$ for every $t$; if $s_\gamma=-1$, then $\gamma$ fixes $(h,0)$. Both are
impossible, since non-trivial deck transformations act freely. So $A_\gamma$ has no fixed point in
$H$ for $\gamma\neq\mathrm{id}$.

\emph{Injectivity.} If $A_\gamma=A_{\gamma'}$, then $\gamma^{-1}\gamma'$ has the form
$(h,t)\mapsto(h,\pm t)$, which fixes $(h,0)$ for all $h$. By freeness, $\gamma^{-1}\gamma'=\mathrm{id}$.

\emph{Proper discontinuity.} Let $K\subset H$ be compact. If $A_\gamma K\cap K\neq\emptyset$, then
$\gamma(K\times[-\tau,\tau])\cap(K\times[-\tau,\tau])\neq\emptyset$. Since $\Gamma$ acts properly
discontinuously on $S$, this happens for finitely many $\gamma$ only, hence, by injectivity, for
finitely many $A_\gamma$.

So $\bar\Gamma$ acts freely and properly discontinuously by isometries on $H$, and $N=H/\bar\Gamma$
is a complete flat manifold, with $H\to N$ a Riemannian covering. The projection $S\to H$ is
equivariant with respect to $\gamma\mapsto A_\gamma$, so it induces a well defined continuous
surjective map $p:M=S/\Gamma\to H/\bar\Gamma=N$. As $M$ is compact, so is $N$.

\medskip

\noindent\emph{Step 5: proof of (3).}

Write $x=[(h,t)]$ and $y=[(h',t')]$. Since $S\to M$ is a Riemannian universal covering with deck
group $\Gamma$,

\[
d_M(x,y)=\min_{\gamma\in\Gamma}\bigl|(h,t)-\gamma(h',t')\bigr|
=\min_{\gamma\in\Gamma}\sqrt{|h-A_\gamma h'|^2+(t-s_\gamma t')^2},
\]

and, since $H\to N$ is a Riemannian covering with deck group $\bar\Gamma$,
$d_N(p(x),p(y))=\min_{\gamma}|h-A_\gamma h'|$. The first formula gives
$d_M(x,y)\ge d_N(p(x),p(y))$. For the upper bound, take $\gamma_0$ realizing
$\min_\gamma|h-A_\gamma h'|$; since $|t-s_{\gamma_0}t'|\le2\tau$,

\[
d_M(x,y)\le\sqrt{|h-A_{\gamma_0}h'|^2+(2\tau)^2}\le d_N(p(x),p(y))+2\tau .
\]

\medskip

\noindent\emph{Step 6: proof of (4).}

We have $\pi(\partial S)=\partial M$, so every component of $\partial M$ is the image of
$P_+=H\times\{\tau\}$ or of $P_-=H\times\{-\tau\}$; thus there are at most two.

(a) Suppose $s\equiv1$, so $\Gamma$ preserves $P_+$ and $P_-$. If $\pi(x)=\pi(y)$ with $x\in P_+$
and $y\in P_-$, then $y=\gamma x$ for some $\gamma\in\Gamma$, which would swap the two hyperplanes,
a contradiction. So $\Sigma^\pm=\pi(P_\pm)$ are disjoint and $\partial M=\Sigma^+\sqcup\Sigma^-$.

If $x,y\in P_+$ and $\pi(x)=\pi(y)$, then $y=\gamma x$ with $\gamma\in\Gamma$; hence
$\pi|_{P_+}:P_+\to\Sigma^+$ is a Riemannian covering with deck group $\Gamma$ acting through
$\gamma\mapsto A_\gamma$, which is free by Step 4. Since $P_+$ is isometric to $H$, we get
$\Sigma^+\cong H/\bar\Gamma=N$ isometrically for the intrinsic metric induced by $g$. In
particular $\Sigma^+$ is compact and connected, and likewise for $\Sigma^-$.

Thus $\partial M$ has exactly two components.

(b) Suppose some $\gamma_1\in\Gamma$ has $s_{\gamma_1}=-1$. Then $s$ is onto $\{\pm1\}$, so
$\Gamma^0=\ker s$ has index $2$, and $\gamma_1P_+=P_-$. Hence $\pi(P_+)=\pi(P_-)=\partial M$,
which is therefore connected. If $x,y\in P_+$ and $\pi(x)=\pi(y)$, then $y=\gamma x$ with
$\gamma P_+\cap P_+\neq\emptyset$, so $\gamma\in\Gamma^0$. Hence $\pi|_{P_+}:P_+\to\partial M$ is a
Riemannian covering with deck group $\Gamma^0$, and $\partial M\cong H/\bar\Gamma^0$
isometrically. By injectivity of $\gamma\mapsto A_\gamma$, $\bar\Gamma^0$ has index $2$ in
$\bar\Gamma$, hence is normal. Therefore $H/\bar\Gamma^0\to H/\bar\Gamma=N$ is a Riemannian double
cover, whose deck group is generated by

\[
\sigma([h])=[A_{\gamma_1}h].
\]

This is well defined because $\bar\Gamma^0$ is normal, and independent of the choice of $\gamma_1$,
since two choices differ by an element of $\Gamma^0$. It is an isometry, and an involution since
$\gamma_1^2\in\Gamma^0$. It is free: if $[A_{\gamma_1}h]=[h]$, then $A_{\gamma_1}h=A_0h$ for some
$A_0\in\bar\Gamma^0$, so $A_0^{-1}A_{\gamma_1}=A_{\gamma}$ with $s_\gamma=-1$ (so $\gamma\neq\mathrm{id}$)
fixes $h$, contradicting Step 4. Finally $(H/\bar\Gamma^0)/\langle\sigma\rangle=H/\bar\Gamma=N$,
and the quotient metric of the free isometric action is the Riemannian metric of $N$.

For the geometric description, let $x_0=\pi(h,\tau)\in\partial M$. The segment
$\{h\}\times[-\tau,\tau]$ meets $\partial S$ only at its endpoints, so it projects to a geodesic of
length $2\tau$ starting orthogonally at $x_0$ and meeting $\partial M$ for the first time at
$\pi(h,-\tau)$; this is the geodesic in the statement, by uniqueness of the normal geodesic. Since
$\gamma_1(h,-\tau)=(A_{\gamma_1}h,\tau)$, we get
$\pi(h,-\tau)=\pi(\gamma_1(h,-\tau))=\pi(A_{\gamma_1}h,\tau)=\sigma(x_0)$.

\end{proof}

\subsection{Orbifolds and flat orbifolds}

\begin{defi}

An \emph{$n$-dimensional orbifold} is a Hausdorff, paracompact topological space
$\mathcal O$ together with a collection $\{(U_i,\widetilde U_i,\Gamma_i,\varphi_i)\}_{i\in I}$,
called an \emph{atlas}, where for each $i$:

\begin{itemize}

\item $U_i\subset\mathcal O$ is open, $\widetilde U_i\subset\mathbb R^n$ is open and connected,
and $\Gamma_i$ is a finite group of diffeomorphisms of $\widetilde U_i$;

\item $\varphi_i:\widetilde U_i\to U_i$ is continuous, $\Gamma_i$-invariant
($\varphi_i\circ\gamma=\varphi_i$ for all $\gamma\in\Gamma_i$), and induces a homeomorphism
$\widetilde U_i/\Gamma_i\to U_i$.

\end{itemize}

Moreover:

\begin{enumerate}

\item the sets $U_i$ cover $\mathcal O$;

\item (compatibility) if $x\in U_i\cap U_j$, there are an index $k$ with $x\in U_k\subset U_i\cap U_j$
and smooth embeddings $\lambda_{ki}:\widetilde U_k\to\widetilde U_i$,
$\lambda_{kj}:\widetilde U_k\to\widetilde U_j$ with
$\varphi_i\circ\lambda_{ki}=\varphi_k=\varphi_j\circ\lambda_{kj}$.

\end{enumerate}

The quadruples $(U_i,\widetilde U_i,\Gamma_i,\varphi_i)$ are called \emph{orbifold charts}.

\end{defi}

Being paracompact, Hausdorff and locally metrizable, $\mathcal O$ is metrizable by Smirnov's metrization theorem.

The group $\Gamma_i$ may contain reflections, so the local models include
in dimension $2$, $\mathbb R^2/\mathbb Z_a$ with $\mathbb Z_a$ generated by a rotation (cone points),
$\mathbb R^2/\langle\text{reflection}\rangle$ (mirrors), and $\mathbb R^2/D_a$ with $D_a$ dihedral of order $2a$
(corner reflectors); see \cite{davis} for orbifolds with reflections. The \emph{underlying space} of $\mathcal O$ is the
topological space $\mathcal O$ itself; it can have boundary even though $\mathcal O$ is not a manifold
with boundary in the orbifold sense. The orbifold is \emph{Riemannian} (resp.\ \emph{flat}) if the
$\Gamma_i$ act by isometries of a Riemannian metric on $\widetilde U_i$ (resp.\ of the Euclidean
metric) and the embeddings $\lambda_{ki}$ are isometric embeddings. A flat orbifold is a length space
with the metric induced by the charts.

\begin{defi}

An orbifold $\mathcal O$ is \emph{good} if it is globally a quotient $\mathcal O=\widetilde{\mathcal O}/\Gamma$,
where $\widetilde{\mathcal O}$ is a manifold and $\Gamma$ is a group acting properly discontinuously by
diffeomorphisms (not necessarily freely); equivalently, $\mathcal O$ admits an orbifold covering by a
manifold. Otherwise it is \emph{bad}. See \cite[Chapter~13]{thurston2}, \cite[Chapter~2]{boilmailport} and \cite[Section~2]{scott} for details about orbifolds; for example, the bad closed $2$-orbifolds are the teardrop, the spindle with cone points of different orders, and their quotients by a reflection (a disc with mirror boundary and one corner reflector, or two corner reflectors of different orders); every other closed $2$-orbifold is good \cite[Section~2]{scott}.
\end{defi}

\begin{prop}\label{prop:good}

A compact connected flat orbifold is good: it is isometric to $\mathbb R^k/\Gamma$, where
$\Gamma\le\mathrm{Isom}(\mathbb R^k)$ is a discrete cocompact subgroup (acting properly
discontinuously, but not necessarily freely), with the quotient metric
$d(\Gamma x,\Gamma y)=\min_{\gamma\in\Gamma}|x-\gamma y|$. Conversely, each such quotient is a compact
flat orbifold.

\end{prop}

\noindent Apply \cite[Chapter~13, Proposition~13.3.2]{thurston2} with $X=\mathbb R^k$ and $G=\mathrm{Isom}(\mathbb R^k)$ (see also \cite[Chapter~13]{ratcliffe}). A compact connected Euclidean orbifold is metrically complete; its developing map identifies its universal orbifold cover with $\mathbb R^k$, and its deck group is discrete and cocompact. This proves the stated quotient description. When $\Gamma$ acts freely we recover the closed flat manifolds discussed above, and Bieberbach's theorem
holds for $\Gamma$ as well (the translation subgroup is a lattice of finite index). In this paper, ``flat orbifold'' always
means a good flat orbifold $\mathbb R^k/\Gamma$ as in Proposition~\ref{prop:good}, and reflections
are allowed.

\medskip

\noindent In dimension $2$, compact flat orbifolds are quotients $\mathbb R^2/\Gamma$, with $\Gamma$ belonging to one of the $17$ affine types of wallpaper (plane crystallographic) groups (see \cite{plesken} for the number of crystallographic groups in dimensions up to six). These are affine types, not individual isometry classes; metric parameters may vary within a type. Table~\ref{tab:flat2orb} lists them using Conway's
orbifold notation \cite{conway}: digits $a$ before a star denote cone points of order $a$, a star $*$ denotes a mirror
boundary circle, digits after a star denote corner reflectors (dihedral groups of order $2a$), and
$\times$ denotes a cross-cap.

\begin{table}[ht]

\centering

\renewcommand{\arraystretch}{1.2}

\begin{tabular}{|c|c|l|p{7cm}|}

\hline
Conway & Wallpaper group & Underlying space & Singularities \\ \hline
$o$ & $p1$ & torus & none \\
$\times\times$ & $pg$ & Klein bottle & none \\
$2222$ & $p2$ & sphere & four cone points of order $2$ \\
$333$ & $p3$ & sphere & three cone points of order $3$ \\
$442$ & $p4$ & sphere & cone points of orders $4,4,2$ \\
$632$ & $p6$ & sphere & cone points of orders $6,3,2$ \\
$22\times$ & $pgg$ & projective plane & two cone points of order $2$ \\
$**$ & $pm$ & annulus & two mirror boundary circles \\
$*\times$ & $cm$ & M\"obius band & one mirror boundary circle \\
$22*$ & $pmg$ & disc & two cone points of order $2$; mirror boundary \\
$*2222$ & $pmm$ & disc & mirror boundary with four corner reflectors of order $2$ \\
$2{*}22$ & $cmm$ & disc & one cone point of order $2$; mirror boundary with two corner reflectors of order $2$ \\
$4{*}2$ & $p4g$ & disc & one cone point of order $4$; mirror boundary with one corner reflector of order $2$ \\
$*442$ & $p4m$ & disc & mirror boundary with corner reflectors of orders $4,4,2$ \\
$3{*}3$ & $p31m$ & disc & one cone point of order $3$; mirror boundary with one corner reflector of order $3$ \\
$*333$ & $p3m1$ & disc & mirror boundary with three corner reflectors of order $3$ \\
$*632$ & $p6m$ & disc & mirror boundary with corner reflectors of orders $6,3,2$ \\ \hline

\end{tabular}
\smallskip

\caption{The $17$ affine types of compact flat $2$-orbifolds.}

\label{tab:flat2orb}

\end{table}

\section{Equivariant Gromov--Hausdorff convergence and proof of Theorem A}\label{theoremA}

\subsection{Equivariant Gromov--Hausdorff convergence}

\begin{defi}

We denote by $\mathcal M$ the set of all isometry classes of pointed
length spaces (in the sense defined in \cite{gromov}) $(X,p)$ such that, for each $D$, the ball $B_p(D,X)$ is relatively
compact. \end{defi}

Let $\mathcal M_{\mathrm{eq}}$ be the set of triples $(X,\Gamma,p)$, where
$(X,p)\in\mathcal M$ and $\Gamma$ is a closed group of isometries of $X$. Put

\begin{equation}
\Gamma(D)=\{\gamma\in\Gamma\mid d(\gamma p,p)<D\}.
\end{equation}

\medskip

\begin{defi}

Let $(X,\Gamma,p),(Y,\Lambda,q)\in\mathcal M_{\mathrm{eq}}$. An $\epsilon$-equivariant pointed Hausdorff approximation stands for a triple
$(f,\varphi,\psi)$ of maps

\[
f:B_p(1/\epsilon,X)\longrightarrow Y,\qquad
\varphi:\Gamma(1/\epsilon)\longrightarrow\Lambda(1/\epsilon)
\]

and

\[
\psi:\Lambda(1/\epsilon)\longrightarrow\Gamma(1/\epsilon)
\]

such that

\begin{enumerate}

 \item $f(p)=q$;

 \item the $\epsilon$-neighborhood of
 $f(B_p(1/\epsilon,X))$ contains $B_q(1/\epsilon,Y)$;

 \item if $x,y\in B_p(1/\epsilon,X)$, then

 \[
 \left|d(f(x),f(y))-d(x,y)\right|<\epsilon;
 \]

 \item if $\gamma\in\Gamma(1/\epsilon)$, $x\in B_p(1/\epsilon,X)$,
 $\gamma x\in B_p(1/\epsilon,X)$, then

 \[
 d\bigl(f(\gamma x),\varphi(\gamma)(f(x))\bigr)<\epsilon;
 \]

 \item if $\mu\in\Lambda(1/\epsilon)$, $x\in B_p(1/\epsilon,X)$,
 $\psi(\mu)(x)\in B_p(1/\epsilon,X)$, then

 \[
 d\bigl(f(\psi(\mu)(x)),\mu f(x)\bigr)<\epsilon.
 \]

\end{enumerate}

\end{defi}

We do not assume either that $f$ is continuous or that $\varphi$ and $\psi$
are homomorphisms.

The \emph{equivariant pointed Hausdorff distance}

\[
d_{\mathrm{ePH}}\bigl((X,\Gamma,p),(Y,\Lambda,q)\bigr)
\]

stands for the infimum of the positive numbers $\epsilon$ such that there exist
$\epsilon$-equivariant pointed Hausdorff approximations from $(X,\Gamma,p)$
to $(Y,\Lambda,q)$ and from $(Y,\Lambda,q)$ to $(X,\Gamma,p)$.
\medskip

By $d_{\mathrm{pH}}$ we denote the pointed Hausdorff distance, which is the
case when the groups are trivial. The notation

\[
\lim_{i\to\infty}(X_i,\Gamma_i,p_i)=(Y,\Lambda,q)
\]

means that

\[
\lim_{i\to\infty}
d_{\mathrm{ePH}}\bigl((X_i,\Gamma_i,p_i),(Y,\Lambda,q)\bigr)=0.
\]

We use equivariant compactness from \cite[Proposition~3.6]{fukayayama}, together with quotient convergence from \cite[Theorem~2-1]{fukaya1}:

\begin{theorem}[Equivariant compactness and quotient convergence]\label{equivariant}

Let $(X_i,p_i)$ be a sequence of pointed proper length spaces converging in the pointed Gromov--Hausdorff sense to a pointed proper length space $(X,p)$. Let $G_i \leq \mathrm{Isom}(X_i)$ be a sequence of closed subgroups.

Then, after passing to a subsequence $k_i$, there exists a closed subgroup $G \leq \mathrm{Isom}(X)$ such that:

\begin{enumerate}

 \item The sequence $(X_{k_i},G_{k_i},p_{k_i})$ converges to $(X,G,p)$ in the equivariant sense.

 \item The pointed quotient spaces converge in the pointed Gromov--Hausdorff sense:

 \[
 \bigl(X_{k_i}/G_{k_i},[p_{k_i}]\bigr)\;\xrightarrow{pGH}\;\bigl(X/G,[p]\bigr).
 \]

\end{enumerate}

\end{theorem}

\begin{rem}

Here $X/G$ carries the quotient metric $d(Gx,Gy)=\inf_{g\in G}d(x,gy)$ ($G$ is closed), and the equivariant approximations in (1) go from $(X_{k_i},G_{k_i},p_{k_i})$ to $(X,G,p)$.

\end{rem}

\begin{lemma}\label{lem:involution_limit}

Let $(X_i,p_i)$ and $(X,p)$ be compact metric spaces with
$\operatorname{diam}X_i\le \Delta$ and $\operatorname{diam}X\le\Delta$, and let
$G_i=\{\mathrm{Id},\iota_i\}\le \mathrm{Isom}(X_i)$ with $\iota_i^2=\mathrm{Id}$
(we allow $\iota_i=\mathrm{Id}$). Assume that
$(X_i,G_i,p_i)\to(X,G,p)$ in the equivariant Gromov--Hausdorff sense, where
$G\le\mathrm{Isom}(X)$ is a closed subgroup. More precisely, let
$(f_i,\varphi_i,\psi_i)$ be $\epsilon_i$-equivariant pointed Hausdorff
approximations from $(X_i,G_i,p_i)$ to $(X,G,p)$, with $\epsilon_i\to0$, so that
$f_i:X_i\to X$, $\varphi_i:G_i\to G$ and $\psi_i:G\to G_i$.

Then, after passing to a subsequence, the following hold.

\begin{enumerate}

\item[\textup{(i)}] There is $\bar\iota\in G$ such that $\varphi_i(\iota_i)\to\bar\iota$
uniformly on $X$, and

\[
\eta_i:=\sup_{x\in X_i} d_X\bigl(f_i(\iota_i x),\bar\iota f_i(x)\bigr)\longrightarrow 0 .
\]

\item[\textup{(ii)}] $\bar\iota^{\,2}=\mathrm{Id}$.

\item[\textup{(iii)}] $G=\{\mathrm{Id},\bar\iota\}$. In particular $|G|\le2$, and
$G$ is trivial if and only if $\bar\iota=\mathrm{Id}$.

\item[\textup{(iv)}] Let $s_i=\sup_{x\in X_i}d(x,\iota_i x)$ and
$m=\max_{y\in X}d(y,\bar\iota y)$. Then $s_i\to m$. In particular,
$\bar\iota=\mathrm{Id}$ if and only if $s_i\to0$.

\end{enumerate}

\end{lemma}

\begin{proof}

Since $\Delta<1/\epsilon_i$ for $i$ large, the balls $B_{p_i}(1/\epsilon_i,X_i)$ and
$B_p(1/\epsilon_i,X)$ are the whole spaces, and $G_i(1/\epsilon_i)=G_i$,
$G(1/\epsilon_i)=G$. So all the maps $f_i,\varphi_i,\psi_i$ are defined on the whole of
their domains, and the properties of an $\epsilon_i$-equivariant approximation hold
globally:

\begin{itemize}

\item[(a)] $f_i(X_i)$ is $\epsilon_i$-dense in $X$ and $|d(f_i(x),f_i(x'))-d(x,x')|<\epsilon_i$;

\item[(b)] $d(f_i(\gamma x),\varphi_i(\gamma)f_i(x))<\epsilon_i$ for $\gamma\in G_i$, $x\in X_i$;

\item[(c)] $d(f_i(\psi_i(g)x),g f_i(x))<\epsilon_i$ for $g\in G$, $x\in X_i$.

\end{itemize}

\medskip

\emph{(i)} Since $X$ is compact, $\mathrm{Isom}(X)$ is compact in the topology of uniform
convergence (Arzel\`a--Ascoli), and $G$ is closed in it, hence compact. As $\varphi_i(\iota_i)\in G$, a subsequence converges uniformly to some
$\bar\iota\in G$. Put $\delta_i=\sup_{y\in X}d(\varphi_i(\iota_i)y,\bar\iota y)\to0$. By (b), for every $x\in X_i$,

\[
d\bigl(f_i(\iota_ix),\bar\iota f_i(x)\bigr)
\le d\bigl(f_i(\iota_ix),\varphi_i(\iota_i)f_i(x)\bigr)+\delta_i<\epsilon_i+\delta_i ,
\]

so $\eta_i\le\epsilon_i+\delta_i\to0$.
\medskip

\emph{(ii)} Let $y\in X$. By (a) choose $x_i\in X_i$ with $d(f_i(x_i),y)<\epsilon_i$. Applying (i) to $x_i$ and to $\iota_ix_i$ gives

\[
d\bigl(f_i(x_i),\bar\iota f_i(\iota_ix_i)\bigr)
=d\bigl(f_i(\iota_i^2x_i),\bar\iota f_i(\iota_ix_i)\bigr)\le\eta_i,
\qquad
d\bigl(f_i(\iota_ix_i),\bar\iota f_i(x_i)\bigr)\le\eta_i .
\]

Since $\bar\iota$ is an isometry, the second inequality gives $d\bigl(\bar\iota f_i(\iota_ix_i),\bar\iota^{\,2}f_i(x_i)\bigr)\le\eta_i$. Hence 

$$d\bigl(f_i(x_i),\bar\iota^{\,2}f_i(x_i)\bigr)\le2\eta_i\to0.$$ 

Because $f_i(x_i)\to y$ and $\bar\iota^{\,2}$ is continuous, $\bar\iota^{\,2}y=y$. As $y$ was arbitrary, $\bar\iota^{\,2}=\mathrm{Id}$.
\medskip

\emph{(iii)} We already know $\bar\iota\in G$. Conversely, let $g\in G$. Then $\psi_i(g)\in\{\mathrm{Id},\iota_i\}$, so after passing to a subsequence (depending on $g$) $\psi_i(g)$ is constant, equal to $\mathrm{Id}$ or to $\iota_i$. Let $y\in X$ and choose $x_i$ as in (ii). If $\psi_i(g)=\mathrm{Id}$, then (c) gives $d(f_i(x_i),gf_i(x_i))<\epsilon_i$, so $gy=y$. Since $y$ is arbitrary, $g=\mathrm{Id}$. If $\psi_i(g)=\iota_i$, then (c) gives $d(f_i(\iota_ix_i),gf_i(x_i))<\epsilon_i$, and with (i)
we get $d(\bar\iota f_i(x_i),gf_i(x_i))\to0$, so $\bar\iota y=gy$ and $g=\bar\iota$. Thus $G\subseteq\{\mathrm{Id},\bar\iota\}$, and equality follows.
\medskip

\emph{(iv)} For $x\in X_i$, by (a) and (i),

\[
\bigl|d(x,\iota_ix)-d\bigl(f_i(x),\bar\iota f_i(x)\bigr)\bigr|
\le \bigl|d(x,\iota_ix)-d(f_i(x),f_i(\iota_ix))\bigr|
+d\bigl(f_i(\iota_ix),\bar\iota f_i(x)\bigr)<\epsilon_i+\eta_i .
\]

Taking the supremum over $x\in X_i$ gives $\bigl|s_i-\sup_{y\in f_i(X_i)}d(y,\bar\iota y)\bigr|\le\epsilon_i+\eta_i$. Since $f_i(X_i)$ is $\epsilon_i$-dense and $y\mapsto d(y,\bar\iota y)$ is $2$-Lipschitz, $\bigl|\sup_{f_i(X_i)}d(y,\bar\iota y)-m\bigr|\le2\epsilon_i$. Hence $|s_i-m|\le3\epsilon_i+\eta_i\to 0$.

\end{proof}

\subsection{Correspondences and the Gromov--Hausdorff distance} We recall the following standard characterization of the Gromov--Hausdorff distance (see \cite{bbi} for details).

\begin{defi} Let $(X,d_X)$ and $(Y,d_Y)$ be metric spaces. A \emph{correspondence} between $X$ and $Y$ is a subset \[ R\subset X\times Y \] such that the two projections \[ R\longrightarrow X \qquad\text{and}\qquad R\longrightarrow Y \] are surjective. Equivalently, for every $x\in X$ there exists $y\in Y$ such that $(x,y)\in R$, and for every $y\in Y$ there exists $x\in X$ such that $(x,y)\in R$. \end{defi}

\begin{defi} Let $R$ be a correspondence between the metric spaces $(X,d_X)$ and $(Y,d_Y)$. The \emph{distortion} of $R$ is defined by \[ \operatorname{dis}(R) = \sup \left\{ \left| d_X(x,x')-d_Y(y,y') \right| : (x,y),(x',y')\in R \right\}. \] \end{defi}

We will use the following characterization of the Gromov--Hausdorff distance.

\begin{theorem} For any metric spaces $X$ and $Y$, \[ d_{GH}(X,Y) = \frac{1}{2} \inf_R \operatorname{dis}(R), \] where the infimum is taken over all correspondences $R$ between $X$ and $Y$. \end{theorem}

\noindent See \cite[Theorem~7.3.25]{bbi}.

\subsection{Proof of Theorem A}

Put $\tau_i=\sup_{x\in M}d_{g_i}(x,\partial M)\le D$ (note $\partial M\ne\emptyset$), and pass to a subsequence along which $\tau_i\to\tau_\infty\in[0,D]$. This is preserved
under further subsequences, so case (i) of the theorem is $\tau_\infty>0$ and case (ii)
is $\tau_\infty=0$. In the arguments below, ``passing to a further subsequence" always
refers to a subsequence of this one.

\noindent For each manifold with geodesic boundary \(M_i=(M,g_i)\), form the double \(\widetilde M_i = M_i\sqcup_{\partial M_i} M_i\). We regard $M_i$ as the first copy. Each $\widetilde M_i$ is closed and connected, since $M$ is compact and connected and $\partial M\ne\emptyset$. Because \(\partial M_i\) is totally geodesic, the doubled metric \(\widetilde g_i\) is flat and smooth (see Section~\ref{flat}: near $\partial M_i$ the metric is a product $dt^2+g_{\partial M_i}$).

Let \(\iota_i:\widetilde M_i\to \widetilde M_i \) denote the canonical involution exchanging the two copies and fixing the glued boundary. The fixed point-set $\operatorname{Fix}(\iota_i)$ of $\iota_i$ is $\operatorname{Fix}(\iota_i)=\partial M_i$, where $\partial M_i$ is identified with the glued hypersurface in $\widetilde M_i$. Moreover, for every $x\in M_i$,

\[
d_{\widetilde g_i}(x,\iota_i(x))
=
2d_{g_i}(x,\partial M_i).
\]

\noindent Here $d_{\widetilde g_i}=d_{g_i}$ on $M_i$: paths in $M_i$ are paths in $\widetilde M_i$, so $d_{\widetilde g_i}\le d_{g_i}$, and the folding map $\widetilde M_i\to M_i$ is 1-Lipschitz and restricts to the identity on $M_i$, so $d_{g_i}\le d_{\widetilde g_i}$. Every path from $x$ to $\iota_ix$ meets $\partial M_i$ at some point $q$ (it passes from one copy to the other), so its length is at least 

\[d_{\widetilde g_i}(x,q)+d_{\widetilde g_i}(q,\iota_ix)=2d_{\widetilde g_i}(x,q)\ge2d_{g_i}(x,\partial M_i),\]

\noindent since $\iota_iq=q$ and $\iota_i$ is an isometry, so $d_{\widetilde g_i}(q,\iota_ix)=d_{\widetilde g_i}(\iota_iq,\iota_ix)=d_{\widetilde g_i}(q,x)$. Equality is attained by the path going from $x$ to a nearest point of $\partial M_i$ and then reflecting by $\iota_i$.
\medskip

\noindent We have that \((\widetilde M_i)\) is a sequence of closed, connected flat manifolds with an isometric involution \(\iota_i\), satisfying $\operatorname{diam}(\widetilde M_i) \leq 2D$. Indeed, let $x,y\in\widetilde M_i$. If they lie in the same copy, then $d_{\widetilde g_i}(x,y)\le d_{g_i}(x,y)\le D$, because paths in a copy are paths in $\widetilde M_i$. If $x\in M_i$ and $y\in\iota_i(M_i)$, let $z\in\partial M_i$ be a point nearest to $x$. Then $d_{\widetilde g_i}(x,z)=d_{g_i}(x,z)\le D$ and, since $\iota_i z=z$ and $\iota_iy\in M_i$, $d_{\widetilde g_i}(z,y)=d_{\widetilde g_i}(z,\iota_iy)=d_{g_i}(z,\iota_iy)\le D$. Hence $d_{\widetilde g_i}(x,y)\le 2D$. The remaining cases follow by symmetry, applying $\iota_i$. Since the metrics are flat, $\mathrm{Ric}\ge0$. By Gromov's precompactness theorem, after passing to a further subsequence, $\widetilde M_i\to O$ in the Gromov--Hausdorff sense for some compact metric space $O$. By \cite{betderpic}, $O$ is a compact flat orbifold with $\dim O=k\le n$ (the hypothesis $\mathrm{Vol}(M,g_i)\to0$ is not used in the proof). Fix $p\in O$ and choose $p_i\in\widetilde M_i$ such that $(\widetilde M_i,p_i)\to(O,p)$ in the pointed Gromov--Hausdorff sense (possible since all spaces are compact). Note that $\operatorname{diam}O\le2D$, as a Gromov--Hausdorff limit of spaces of diameter $\le2D$. The groups $G_i=\{\mathrm{Id},\iota_i\}$ are finite, hence closed in $\mathrm{Isom}(\widetilde M_i)$. Applying Theorem~\ref{equivariant} and passing to a further subsequence, we obtain a closed (hence compact) group $G\le\mathrm{Isom}(O)$ with $(\widetilde M_i,G_i,p_i)\to(O,G,p)$ equivariantly, with $\epsilon_i$-equivariant approximations $(f_i,\varphi_i,\psi_i)$ from $(\widetilde M_i,G_i,p_i)$ to $(O,G,p)$, $\epsilon_i\to0$. By Lemma~\ref{lem:involution_limit} (with $\Delta=2D$), after a further subsequence there is $\bar\iota\in G$ with $\bar\iota^2=\mathrm{Id}$ and $G=\{\mathrm{Id},\bar\iota\}$, so $|G|\le2$. We fix this subsequence, the approximations and $\bar\iota$ from now on. We need the following lemma.

\begin{lemma}\label{lem:trivial_action}

In the situation above, with $\bar\iota$ as in Lemma~\ref{lem:involution_limit},

\[
G=\{\mathrm{Id}\}\iff\bar\iota=\mathrm{Id}\iff\tau_\infty=0 .
\]

\end{lemma}

\begin{proof}

We first compute $s_i:=\sup_{x\in\widetilde M_i}d_{\widetilde g_i}(x,\iota_ix)$. For $x\in M_i$ (the first copy) we have $d_{\widetilde g_i}(x,\iota_ix)=2d_{g_i}(x,\partial M_i)$,
so the supremum over the first copy equals $2\tau_i$. If $y$ lies in the second copy, then
$x:=\iota_iy\in M_i$ and $d(y,\iota_iy)=d(\iota_ix,x)$, so the supremum over the second copy
is also $2\tau_i$. (Points of $\partial M_i$ belong to both copies and have
$d(y,\iota_iy)=0$, so there is no inconsistency.) Hence

\[
s_i=2\tau_i .
\]

By construction $\tau_i\to\tau_\infty$ along our subsequence, so
$s_i\to2\tau_\infty$. On the other hand, Lemma~\ref{lem:involution_limit}(iv) gives
$s_i\to m:=\max_{y\in O}d(y,\bar\iota y)$. Therefore

\[
m=2\tau_\infty .
\]

Now $m=0$ if and only if $\bar\iota=\mathrm{Id}$, so $\bar\iota=\mathrm{Id}$ if and only if
$\tau_\infty=0$. Finally, by Lemma~\ref{lem:involution_limit}(iii) we have
$G=\{\mathrm{Id},\bar\iota\}$, hence $G=\{\mathrm{Id}\}$ if and only if $\bar\iota=\mathrm{Id}$.

This proves both equivalences.

\end{proof}

\begin{enumerate}

\medskip

\item[\textup{(i)}] \emph{Suppose that $\tau_\infty>0$.} By Lemma~\ref{lem:trivial_action}, $\bar\iota\neq\mathrm{Id}$ and $G=\{\mathrm{Id},\bar\iota\}\cong\mathbb Z_2$. By Lemma~\ref{lem:involution_limit}(ii), $\bar\iota^{2}=\mathrm{Id}$. Thus $\bar\iota$ is a
non-trivial isometric involution of the compact flat orbifold $O$. Let $\pi_i\colon\widetilde M_i\to\widetilde M_i/G_i$ be the quotient map. We first check that $\pi_i|_{M_i}\colon M_i\to\widetilde M_i/G_i$ is an isometry, where $\widetilde M_i/G_i$ carries the quotient metric $d(\pi_i x,\pi_i y)=\min\{d_{\widetilde g_i}(x,y),d_{\widetilde g_i}(x,\iota_iy)\}$. It is bijective, since every $\langle\iota_i\rangle$-orbit meets $M_i$ in exactly one point (the orbit of an interior point $x\in M_i\setminus\partial M_i$ is $\{x,\iota_ix\}$ with $\iota_ix\notin M_i$, and points of $\partial M_i$ are fixed). For $x,y\in M_i$ we have $d_{\widetilde g_i}(x,y)=d_{g_i}(x,y)$. Moreover, the folding map $r:\widetilde M_i\to M_i$ is 1-Lipschitz with $r(x)=x$ and $r(\iota_iy)=y$, so $d_{\widetilde g_i}(x,\iota_iy)\ge d_{g_i}(x,y)$. Hence the minimum equals $d_{g_i}(x,y)$. Thus $(M_i,g_i)\cong(\widetilde M_i/G_i,\widetilde g_i/G_i)$ as metric spaces. Hence, by the quotient convergence in Theorem~\ref{equivariant}(2), along our fixed subsequence

$$
(M_i,g_i)
\cong
(\widetilde M_i/G_i,\widetilde g_i/G_i)
\xrightarrow{GH}
O/G
=
O/\langle\bar\iota\rangle.
$$

\noindent Consequently the Gromov--Hausdorff limit of $(M_i,g_i)$ is $X=O/\langle\bar\iota\rangle$. It remains to see that $X$ is a compact flat orbifold (a good orbifold $\mathbb R^k/\Gamma^{\sharp}$ with $\Gamma^{\sharp}\le\mathrm{Isom}(\mathbb R^k)$ discrete and cocompact, possibly containing reflections). Since $O$ is a compact flat orbifold, it is good: $O=\mathbb R^k/\Gamma$ with $\Gamma\le\mathrm{Isom}(\mathbb R^k)$ discrete and cocompact (\cite[Chapter~13, Proposition~13.3.2]{thurston2}). A metric isometry of a Riemannian orbifold is an orbifold isometry; the local lifting result and the universal orbifold-cover property give an isometric lift $\tilde\iota$ of $\bar\iota$ to $\mathbb R^k$ (\cite[Lemma~2.2, Remark~2.3 and Theorem~2.9]{lange}). Conjugating deck transformations by this lift gives $\tilde\iota\Gamma\tilde\iota^{-1}=\Gamma$. Since $\bar\iota^2=\mathrm{Id}$, for every $x$ we have $\tilde\iota^2x\in\Gamma x$. The closed sets $F_\gamma=\{x:\tilde\iota^2x=\gamma x\}$ cover $\mathbb R^k$ and $\Gamma$ is countable, so by Baire some $F_\gamma$ has interior; as $\tilde\iota^2$ and $\gamma$ are affine, $\tilde\iota^2=\gamma\in\Gamma$. Hence $\Gamma^{\sharp}=\Gamma\cup\tilde\iota\Gamma$ is a group since $\tilde\iota$ normalizes $\Gamma$ and $\tilde\iota^2\in\Gamma$, so $\Gamma\cup\tilde\iota\Gamma$ is closed under products and inverses. It contains $\Gamma$ with index $\le2$, so it is discrete and cocompact and $X=O/\langle\bar\iota\rangle=\mathbb R^k/\Gamma^{\sharp}$ as metric spaces, since both carry the quotient metric $d(\Gamma^{\sharp}x,\Gamma^{\sharp}y)=\min_{\gamma^{\sharp}\in\Gamma^{\sharp}}|x-\gamma^{\sharp}y|$. Therefore $X=\mathbb R^k/\Gamma^{\sharp}$ is a compact flat orbifold.
\medskip

\item[\textup{(ii)}] \emph{Suppose that $\tau_\infty=0$.} By Lemma~\ref{lem:trivial_action}, $\bar\iota=\mathrm{Id}$ and $G=\{\mathrm{Id}\}$, and $\tau_i\to0$ along our fixed subsequence. Recall that $d_{\widetilde g_i}=d_{g_i}$ on $M_i$, so $M_i\subset\widetilde M_i$ is an isometric embedding. For every $x\in\widetilde M_i$ there is $q\in\partial M_i\subset M_i$ with $d_{\widetilde g_i}(x,q)=d_{\widetilde g_i}(x,\partial M_i)\le\tau_i$ (for $x$ in the second copy apply $\iota_i$: since $\iota_i$ is an isometry fixing $\partial M_i$, $d_{\widetilde g_i}(x,\partial M_i)=d_{\widetilde g_i}(\iota_ix,\partial M_i)\le\tau_i$, because $\iota_ix\in M_i$), so $M_i$ is $\tau_i$-dense in $\widetilde M_i$. Since $M_i\subset\widetilde M_i$ is an isometric embedding, the Hausdorff distance between $M_i$ and $\widetilde M_i$, viewed as subsets of $\widetilde M_i$, is at most $\tau_i$, and therefore $d_{GH}(M_i,\widetilde M_i)\le d_H(M_i,\widetilde M_i)\le\tau_i\to0$ (see \cite{bbi}, Chapter~7). Since $\widetilde M_i\to O$, we get $(M_i,g_i)\xrightarrow{GH}O$. In particular $X=O$, a compact flat orbifold. Apply Theorem~\ref{theorem:slab} to each $(M,g_i)$ (here $n\ge2$, and $\tau_i>0$ for every $i$ by the theorem, so it applies along our whole subsequence). Let $S_i=H\times[-\tau_i,\tau_i]$, with $H=\mathbb R^{n-1}$, be the universal cover of $(M,g_i)$, let $\bar\Gamma_i$ be the group of Theorem~\ref{theorem:slab}(2) associated with $g_i$, let $N_i=H/\bar\Gamma_i$, and let $\varpi_i\colon M\to N_i$ be the map $p$ of Theorem~\ref{theorem:slab}(3) for $g_i$. By Theorem~\ref{theorem:slab}(3), the correspondence $R_i=\{(x,\varpi_i(x)):x\in M\}$ between $M_i$ and $N_i$ (it is a correspondence because $\varpi_i$ is onto)
has distortion $\operatorname{dis}(R_i)\le2\tau_i$. By the characterization of the Gromov--Hausdorff distance via correspondences, $d_{GH}(M_i,N_i)\le\tau_i\to0$, hence $N_i\xrightarrow{GH}X$ as well.

\medskip

\noindent\emph{$\partial M$ disconnected.} The number of connected components of $\partial M$ depends only on the manifold $M$, not on the metric $g_i$, so the cases ``connected'' and ``disconnected'' are the same for every $i$. Suppose $\partial M$ is disconnected. By Theorem~\ref{theorem:slab}(4), $\partial M$ has exactly two components $\Sigma_i^\pm$ (this is the case $s\equiv1$ of the theorem, applied to $g_i$), and each of them, with its intrinsic metric $g_i|_{\Sigma_i^\pm}$, is isometric to $N_i$ by Theorem~\ref{theorem:slab}(4a). Hence, by the estimate above,

\[
d_{GH}\bigl((\Sigma_i^\pm,g_i|_{\Sigma_i^\pm}),M_i\bigr)=d_{GH}(N_i,M_i)\le\tau_i\to0 .
\]

\noindent Since $(M_i,g_i)\xrightarrow{GH}X$, we conclude that

\[
X=\lim_{GH}\bigl(\Sigma_i^\pm,g_i|_{\Sigma_i^\pm}\bigr)
\]

for either component.

\medskip

\noindent\emph{$\partial M$ connected.} Put $B_i=(\partial M,g_i|_{\partial M})$, with its intrinsic metric. By Theorem~\ref{theorem:slab}(4b), $B_i$ carries a free isometric involution $\sigma_i$, and its metric quotient $B_i/\langle\sigma_i\rangle$ is isometric to $N_i$. Consequently,
\[
X=\lim_{GH}\bigl(B_i/\langle\sigma_i\rangle\bigr).
\]

\end{enumerate}

\begin{rem}\label{rem:boundary_limits}
In the degenerate case $\tau_i\to0$, the boundary equipped with the restriction of the distance of the double is $\tau_i$-dense in $\widetilde M_i$ and therefore has the same Gromov--Hausdorff limit. If the boundary is disconnected, each component with its intrinsic metric also has this limit, by Theorem~A(ii)(a). If the boundary is connected, the intrinsic metric quotients $B_i/\langle\sigma_i\rangle$ converge to $X$, whereas the intrinsic boundary spaces $B_i$ need not: the M\"obius-band example in the introduction has boundary circles of length~$2$ and limit a circle of length~$1$.

More generally, along the equivariant convergence used above, the boundary images converge as subsets of the limiting double to $F=\operatorname{Fix}(\bar\iota)$. To see this, take approximations $f_i:\widetilde M_i\to O$ with distortion and density error at most $\epsilon_i\to0$, and put $\eta_i=\sup_x d_O(f_i(\iota_i x),\bar\iota f_i(x))\to0$. The displacement identity, valid on both copies of $M_i$, gives
\[
\left|d_O(f_i(x),\bar\iota f_i(x))-2d_{\widetilde g_i}(x,\partial M_i)\right|\le\epsilon_i+\eta_i.
\]
For $b_i\in\partial M_i$, the displacement of $f_i(b_i)$ under $\bar\iota$ is at most $\eta_i$. Compactness of $O$ implies that $F$ is nonempty and that $\sup_{b\in\partial M_i}d_O(f_i(b),F)\to0$. Conversely, if $y\in F$, choose $x_i$ with $d_O(f_i(x_i),y)\le\epsilon_i$. Then $d_{\widetilde g_i}(x_i,\partial M_i)\le(3\epsilon_i+\eta_i)/2$. Choosing a nearest boundary point $b_i$ yields $d_O(f_i(b_i),y)\le(7\epsilon_i+\eta_i)/2$. Hence
\[
d_H^O\bigl(\overline{f_i(\partial M_i)},F\bigr)\longrightarrow0.
\]
Statements below identifying the limit of the boundary with $F$ use this ambient meaning, not the intrinsic boundary distance.
\end{rem}

\section{Flat manifolds with geodesic boundary in dimension 3}\label{klein}

Most of the results of this section are obtained in \cite{elaf}. By Theorem~\ref{theorem:slab}, the universal cover of a compact flat $n$-manifold $M$ with non-empty totally geodesic boundary is isometric to $[-d,d]\times\mathbb R^{n-1}$, where $d=\tau=\sup_{x\in M}d(x,\partial M)>0$.

\bigskip

Now, we study the manifolds themselves. Let $I=[-d,d]$ and $\Gamma=\pi_1(M)$ be the group of deck transformations of the universal covering

\[
I\times\mathbb R^{n-1}\longrightarrow M.
\]

Viewed as a subgroup of $\operatorname{Isom}(\mathbb R^n)$, the group $\Gamma$ preserves the set

\[
\{-d\}\times\mathbb R^{n-1}
\cup
\{d\}\times\mathbb R^{n-1}.
\]

Moreover, the hyperplane $P=\{0\}\times\mathbb R^{n-1}$ is $\Gamma$-invariant, and the induced action of $\Gamma$ on $P$ is
cocompact. Two cases are possible.

\begin{enumerate}

\item Each connected component of the boundary of $I\times\RR^{n-1}$ is $\Gamma$-invariant. Then, any parallel hyperplane
to these components is also $\Gamma$-invariant. In this situation, $\Gamma$ is isomorphic to its restriction to
$\RR^{n-1}$, which is a $(n-1)$-dimensional Bieberbach group that we denote $\Gamma'$. Then $M$ is isometric to $I\times(\RR^{n-1}/\Gamma')$. Its boundary has two connected components and can be identified with $\RR^{n-1}/\Gamma'\sqcup \RR^{n-1}/\Gamma'$.

\item There are elements of $\Gamma$ which exchange the components. Then the subgroup $\Gamma_0$ of $\Gamma$ which preserves these components is a normal subgroup of index $2$. We have the following property.

\end{enumerate}

\begin{lemma}

The restriction map $\gamma\mapsto\gamma_{\mid P}$ is injective.

\end{lemma}

\begin{proof}

Suppose that $\gamma_{\mid P}$ is the identity map. If
$\gamma\in\Gamma_0$, then $\gamma$ preserves each boundary component,
and fixes $P$ pointwise. The only isometries of $\mathbb R^n$ fixing a hyperplane pointwise are
the identity and the reflection in $P$, and the latter exchanges the two sides of $P$, so it does
not preserve the boundary components. Hence $\gamma=\mathrm{id}$.

If $\gamma\in\Gamma\setminus\Gamma_0$, then $\gamma$ exchanges the two
boundary components. Since $\gamma_{\mid P}=id_P$, it follows that
$\gamma$ is the reflection with respect to $P$. This is impossible,
since $\Gamma$ acts freely.

\end{proof}

The restriction map therefore identifies $\Gamma$ and $\Gamma_0$
with their respective images $\Gamma'$ and $\Gamma_0'$ in
$\operatorname{Isom}(P)$. The image $\Gamma'$ is discrete and
cocompact in $\operatorname{Isom}(P)$ (see Theorem~\ref{theorem:slab}(2)), and it acts freely: indeed, if
$\gamma'=\gamma_{\mid P}$ fixed a point $y\in P$, then $\gamma$
would fix $(0,y)$, contradicting the freeness of the action of
$\Gamma$. Thus $\Gamma'$ and $\Gamma_0'$ are Bieberbach groups of
dimension $n-1$.

Denoting by $(x,y)$ the elements of $\RR\times\RR^{n-1}$,
the actions of $\Gamma $ on $\RR^n$ and of $\Gamma^\prime$ on $\RR^{n-1}$ are related
as follows 

$$ \gamma\cdot (x,y)=\left\{\begin{array}{rl} (x,\gamma'\cdot y)&\mathrm{if}\ \gamma\in\Gamma_0\\
(-x,\gamma'\cdot y)&\mathrm{if}\ \gamma\notin\Gamma_0\end{array}\right.$$

The boundary of $M$ is connected and isometric to 
$\RR^{n-1}/\Gamma_0'$, and $M$ admits a Riemannian covering of order $2$
isometric to $I\times(\RR^{n-1}/\Gamma_0')$. This situation is exactly what happens in dimension $2$
with the M\"obius band.

\begin{rem}

There is a strong contrast with the hyperbolic case. By Theorem~\ref{theorem:slab}(4), a compact
flat manifold with non-empty totally geodesic boundary has at most two boundary components.

In dimension $2$, a compact orientable hyperbolic surface of genus $g$ with $k\ge1$ geodesic
boundary components exists if and only if $2-2g-k<0$ (see \cite{buser1}). Hence the number of boundary components
is not bounded: for every $k\ge1$ there are such surfaces for $g\ge1$ (for instance, the one-holed torus has
connected geodesic boundary), whereas for $g=0$ one needs $k\ge3$, as for a pair of pants.

For hyperbolic $3$-manifolds with geodesic boundary, see \cite{fp}.

\end{rem}

\subsection{The three dimensional case}\label{classificationdim3}

In this section and the next we use the notation $t_w$ for the translation by $w$,
$\sigma_{P,w}$ for the reflection in a plane $P$ composed with the translation by a vector $w$
parallel to $P$, and $\sigma_{D,w}$ for the half-turn about a line $D$ composed with the
translation by a vector $w$ parallel to $D$.

In the case 1, i.e., if each connected component of the boundary of $I\times\RR^{2}$ is $\Gamma$-invariant, then $M$ is the product of a flat torus or a flat Klein bottle by a segment.

In the case 2, i.e., if the subgroup of $\Gamma$ which preserves the connected components of the boundary of $I\times\RR^{2}$ is a subgroup $\Gamma_0$ of index $2$, then $\Gamma'$ is either a lattice of the Euclidean plane or the fundamental group of a Klein bottle.

\medskip

\begin{itemize}

\item If $\Gamma'$ is a lattice, then $\Gamma_0'$ is an index-$2$ sublattice of $\Gamma'$. We may choose a basis $(v,v')$ of $\Gamma_0'$ such that $\Gamma'=\mathbb Zv+\mathbb Z\frac{v'}{2}$ (an index-$2$ superlattice of $\Gamma_0'$ is obtained by adding $\frac v2$, $\frac{v'}2$ or $\frac{v+v'}2$, and in each case the basis can be adapted). The manifold is the quotient of $[-d,d]\times\RR^2$ by the subgroup of $\mathrm{Isom}(\RR^3)$ generated by $t_v$ and the glide reflection $\sigma_{P,\frac{v^\prime}{2}}$, where $P$ is the plane $x=0$ with $v$ and $v'$ parallel to $P$.

This manifold is not orientable; its boundary is the flat torus
$\RR^2/(\ZZ v+\ZZ v^\prime)$, whereas the image of $P$ in the quotient is the flat torus

\[
\mathbb R^2/\left(\mathbb Z v+\mathbb Z\frac{v'}2\right).
\]

\noindent The boundary and the soul are both tori. We denote it by $\TT_{v,v^\prime,d}$. Note that such a manifold exists in any dimension.

\medskip

\item If $\Gamma'$ is the Klein-bottle group of $K_{a,b}$, write $\alpha(y,z)=(y+a/2,-z)$ and $\beta(y,z)=(y,z+b)$. The relation $\alpha\beta\alpha^{-1}=\beta^{-1}$ shows that there are three nontrivial homomorphisms $\Gamma'\to\mathbb Z_2$. Their kernels are $\langle\alpha^2,\beta\rangle$, $\langle\alpha,\beta^2\rangle$ and $\langle\alpha\beta,\beta^2\rangle$. The translation $c(y,z)=(y,z+b/2)$ normalizes $\Gamma'$ and satisfies $c\alpha c^{-1}=\beta\alpha$ and $c\beta c^{-1}=\beta$. It conjugates the second kernel to the third, so there are two types of index-$2$ subgroups up to this isometric equivalence.

\medskip

\begin{itemize}

 \item First, $\Gamma_0'$ is the translation lattice of $\Gamma'$ (the subgroup $\langle\alpha^2,\beta\rangle$, where $\alpha$ is the glide reflection and $\beta$ the translation). The manifold is then the quotient of $[-d,d]\times\RR^2$ by the subgroup of $\mathrm{Isom}(\RR^3)$
generated by $t_{(0,0,b)}: (x,y,z)\mapsto (x,y,z+b)$ and $\sigma_{D,(0,\frac{a}{2},0)}:(x,y,z)\mapsto (-x,y+\frac{a}{2},-z)$ where $D$ is the straight line $x=z=0$. This manifold is orientable, and its boundary is the rectangular torus $T_{a,b}=\mathbb R^2/(a\mathbb Z\oplus b\mathbb Z)$. The soul is the Klein bottle $K_{a,b}$. We denote it by $\TK_{a,b,d}$.

\medskip 

\item

Finally, take $\Gamma_0'=\langle\alpha,\beta^2\rangle$, generated on the central plane by $(y,z)\mapsto(y+a/2,-z)$ and $(y,z)\mapsto(y,z+2b)$.

The manifold is then the quotient of $[-d,d]\times\RR^2$ by the subgroup of $\mathrm{Isom}(\RR^3)$
generated by
$\sigma_{P,(0,0,b)}: (x,y,z)\mapsto (-x,y,z+b)$ and $\sigma_{P',(0,\frac{a}{2},0)}\colon(x,y,z)\mapsto (x,y+\frac{a}{2},-z)$
where $P$ and $P'$ are respectively the planes $x=0$ and $z=0$. 

This manifold is not orientable; the boundary and the soul are both Klein bottles.

We denote it by $\KK_{a,b,d} $.

\end{itemize}

\end{itemize}

Therefore, we have three types of compact flat manifolds with connected geodesic boundary. They can be viewed as three-dimensional analogues of the M\"obius band. Following \cite{elaf}, we call them \emph{solid M\"obius bands}.

The following table summarizes this discussion.

\begin{table}[htbp]\renewcommand{\arraystretch}{1.2}

\centering

\begin{tabular}{|c|c|c|c|c|}

\hline
Manifold & Generators of $\Gamma$ & Orientable & Boundary & Soul $\mathbb R^2/\Gamma'$ \\ \hline
$T^2\times[-d,d]$ & $t_v,t_{v'}$ & yes & $T^2\sqcup T^2$ & $T^2$ \\
$K\times[-d,d]$ & glide, translation & no & $K\sqcup K$ & $K$ \\ \hline
$\TT_{v,v',d}$ & $t_v,\ \sigma_{P,v'/2}$ & no & $\mathbb R^2/(\mathbb Zv+\mathbb Zv')$ & $\mathbb R^2/(\mathbb Zv+\mathbb Z\frac{v'}2)$ \\
$\TK_{a,b,d}$ & $t_{(0,0,b)},\ \sigma_{D,(0,\frac a2,0)}$ & yes & $T_{a,b}$ & $K_{a,b}$ \\
$\KK_{a,b,d}$ & $\sigma_{P,(0,0,b)},\ \sigma_{P',(0,\frac a2,0)}$ & no & $K_{a,2b}$ & $K_{a,b}$ \\ \hline

\end{tabular}
\smallskip

\caption{Compact flat $3$-manifolds with totally geodesic boundary.}

\end{table}

These manifolds will be useful later when discussing doubles. In particular, doubling any of them along its totally geodesic boundary produces a closed flat $3$-manifold; see Section~\ref{theoB}.

\section{Collapse of flat $3$-manifolds with geodesic boundary}\label{theoB}

There are ten diffeomorphism types of closed flat $3$-manifolds, six orientable and four
non-orientable. The doubles of the two products
$T^2\times[-d,d]$ and $K\times[-d,d]$ of Section~\ref{classificationdim3} are $T^3$ and
$K\times S^1$. We show that the doubles of the three manifolds with \emph{connected} geodesic
boundary ($\TT$, $\TK$, $\KK$) are the flat manifolds $B_1$, $G_2$, $B_3$ in the notation of \cite[Section~3.5]{wolfspaces}. Note that $B_1$ is
diffeomorphic to $K\times S^1$, since its group is
$\langle\sigma_{P,a_1/2},t_{a_3}\rangle\times\langle t_{a_2}\rangle$.

Altogether, exactly four of the ten closed flat $3$-manifolds arise as doubles of compact flat
$3$-manifolds with totally geodesic boundary: $T^3$, $K\times S^1$, $G_2$ and $B_3$. These four are
pairwise non-diffeomorphic: $T^3$ has trivial holonomy, $B_3$ has holonomy $\mathbb Z_2\times\mathbb Z_2$,
and $K\times S^1$ and $G_2$ both have holonomy $\mathbb Z_2$ but $K\times S^1$ is non-orientable while $G_2$
is orientable.

\subsection{Gluing along the boundary}\label{recollement}

It is well known that the Klein bottle is obtained by gluing two isometric M\"obius bands along their
boundary. A similar phenomenon occurs in higher dimensions. We first give a general criterion.

\begin{lemma}\label{lem:double}

Let $L>0$, $a_3=(L,0,0)$, and let $\widetilde\Gamma\le\mathrm{Isom}(\mathbb R^3)$ be a discrete group
acting freely, with $\mathbb R^3/\widetilde\Gamma$ compact, containing $t_{a_3}$, such that every element has the form

\[
g(x,p)=(\epsilon x+kL,\;g'p),\qquad \epsilon=\pm1,\ k\in\mathbb Z,\ g'\in\mathrm{Isom}(\mathbb R^2),
\]

and such that $\widetilde\Gamma$ contains an element with $\epsilon=-1$ and $x$-part $x\mapsto L-x$.

Put $S=[\frac L4,\frac{3L}4]\times\mathbb R^2$, let $\Gamma_S\le\widetilde\Gamma$ be the stabilizer of $S$
(its elements preserve or exchange the two boundary planes of $S$), and let $\Gamma'_S$ be the group
of restrictions of the elements of $\Gamma_S$ to the central plane $\{x=\frac L2\}\cong\mathbb R^2$.
Let $\Gamma'_{S,0}$ be the restrictions of the elements preserving each boundary plane of $S$.

Then $\Gamma'_{S,0}$ has index $2$ in $\Gamma'_S$, and $S/\Gamma_S$ is the compact flat $3$-manifold with
connected totally geodesic boundary and half-width $d=L/4$ whose universal cover is the slab
$[-d,d]\times\mathbb R^2$ and whose deck group acts on the soul plane through $\Gamma'_S$, with
$\Gamma'_{S,0}$ the subgroup preserving the two boundary planes. By the classification of
Section~\ref{classificationdim3}, it is $\TT_{v,v',d}$ if $\Gamma'_S$ is a lattice, $\TK_{a,b,d}$ if
$\Gamma'_S$ is a Klein-bottle group and $\Gamma'_{S,0}$ is its translation lattice, and
$\KK_{a,b,d}$ if $\Gamma'_S$ is a Klein-bottle group and $\Gamma'_{S,0}$ is another index-$2$ subgroup
(up to conjugation by a translation of $\mathbb R^2$). Moreover,
\[
\mathbb R^3/\widetilde\Gamma\;\cong\;\mathrm{Double}\bigl(S/\Gamma_S\bigr).
\]

\end{lemma}

\begin{proof}
Let $\rho(x,p)=(\frac L2-x,p)$, the reflection in $\{x=\frac L4\}$. For $g$ as above,
$\rho g\rho^{-1}(x,p)=(\epsilon x+k'L,\,g'p)$ with $k'=-k$ if $\epsilon=1$ and $k'=1-k$ if
$\epsilon=-1$, so $\rho g\rho^{-1}=t_{a_3}^{m}g\in\widetilde\Gamma$ with $m=-2k$ or $m=1-2k$.
Hence $\rho$ normalizes $\widetilde\Gamma$ and induces an isometric involution $\iota$ of
$B=\mathbb R^3/\widetilde\Gamma$.

\emph{Tiling.} The $\widetilde\Gamma$-orbit of $\{x=\frac L4\}$ is the set of planes
$\{x=\frac L4+j\frac L2\}$, $j\in\mathbb Z$ (it contains $\{x=\frac{3L}4\}$ and it is invariant under $t_{a_3}^{\pm1}$). They cut
$\mathbb R^3$ into slabs of width $\frac L2$, alternately translates of $S$ and of $\rho(S)=[-\frac L4,\frac L4]\times\mathbb R^2$ by multiples of $L$.
Every element of $\widetilde\Gamma$ permutes the slabs and sends the centre $\frac L2$ of $S$ to a point of $\frac L2+L\mathbb Z$, whereas the
centre of $\rho(S)$ is $0\notin\frac L2+L\mathbb Z$. So no element maps $S$ to a translate of $\rho(S)$.
The planes of the orbit project to one totally geodesic hypersurface $T\subset B$, because
$\{x=\frac L4\}$ and $\{x=\frac{3L}4\}$ are exchanged by an element with $x$-part $x\mapsto L-x$. Cutting $B$ along $T$
therefore gives exactly two components, the images of $S$ and of $\rho(S)$. Moreover $\rho$ fixes $\{x=\frac L4\}$ pointwise and
acts on $\{x=\frac{3L}4\}$ as $t_{a_3}^{-1}$, so $\iota$ fixes $T$ pointwise; it exchanges the two components.

\emph{The first component.} An element of $\widetilde\Gamma$ sending an interior point of $S$ into $S$ maps the tile $S$ to itself, and
$x\mapsto\epsilon x+kL$ preserves $[\frac L4,\frac{3L}4]$ only for $x\mapsto x$ and $x\mapsto L-x$. Hence $\Gamma_S$ acts by
$(x,p)\mapsto(x,g'p)$ or $(x,p)\mapsto(L-x,g'p)$, and the first component is $S/\Gamma_S$, which is compact because $B$ is.
Its boundary is the image $T$ of the two boundary planes, connected since $\Gamma_S$ contains an element exchanging them.

\emph{The group $\Gamma'_S$.} The restriction $\Gamma_S\to\Gamma'_S$ is injective: a nontrivial element restricting to the identity on $\{x=\frac L2\}$
would be the reflection in this plane, which has fixed points. The homomorphism $\Gamma_S\to\mathbb Z_2$ recording whether the boundary planes are exchanged
is onto, with kernel $\Gamma_{S,0}$, so $\Gamma'_{S,0}$ has index $2$ in $\Gamma'_S$. The group $\Gamma'_S$ acts freely on $\mathbb R^2$: if $g'p=p$, then $g$ fixes
$(\frac L2,p)$ in both cases, a contradiction. It is cocompact because $S/\Gamma_S$ is compact, and discrete because $\Gamma_S$ acts properly discontinuously on $S$.
So $\Gamma'_S$ is a Bieberbach group of dimension $2$. After the translation $x\mapsto x-\frac L2$, the slab $S$ becomes $[-d,d]\times\mathbb R^2$ with $d=\frac L4$
and $x\mapsto L-x$ becomes $x\mapsto -x$; thus $S/\Gamma_S$ is the manifold of Section~\ref{classificationdim3} with half-width $d=\frac L4$ and the pair
$(\Gamma'_S,\Gamma'_{S,0})$, and the classification there gives $\TT$, $\TK$ or $\KK$ as stated.

\emph{The double.} The second component is $\iota(S/\Gamma_S)$, and $\iota$ is an isometry fixing $T=\partial(S/\Gamma_S)$ pointwise. So $B$ is
$S/\Gamma_S$ glued to a copy of itself along its boundary by the identity, i.e.\ $B\cong\mathrm{Double}(S/\Gamma_S)$.
\end{proof}

We now apply the lemma to three groups.

\begin{itemize}

\item \textbf{The double of $\TT_{v,v',d}$.} Let $P=\{x=0\}$, let $a_1,a_2\in P$ be linearly independent and
$a_3=(L,0,0)$. Put

\[
\widetilde\Gamma_{B_1}=\langle\sigma_{P,a_1/2},\,t_{a_2},\,t_{a_3}\rangle,\qquad B_1=\mathbb R^3/\widetilde\Gamma_{B_1},
\]

\noindent where $\sigma_{P,a_1/2}(x,p)=(-x,p+\frac{a_1}2)$ and $\sigma_{P,a_1/2}^{\,2}=t_{a_1}$. The group satisfies the
hypotheses of Lemma~\ref{lem:double}: the element $g=t_{a_3}\circ\sigma_{P,a_1/2}$ has $x$-part
$x\mapsto L-x$. The group $\Gamma_S$ is generated by $t_{a_1}$, $t_{a_2}$ (which preserve the boundary planes
of $S$) and $g$ (which exchanges them). On the central plane $\{x=\frac L2\}$, $g$ acts as the translation by
$\frac{a_1}2$, so

\[
\Gamma'_S=\mathbb Za_2+\mathbb Z\tfrac{a_1}2,\qquad \Gamma'_{S,0}=\mathbb Za_2+\mathbb Za_1 .
\]

These are the groups defining $\TT_{v,v',d}$ with $v=a_2$, $v'=a_1$, and $d=L/4$. Hence

\[
B_1\cong\mathrm{Double}\bigl(\TT_{a_2,a_1,L/4}\bigr).
\]

\item \textbf{The double of $\TK_{a,b,d}$.} Let $D=\{x=z=0\}$, $a_1=(0,0,b)$, $a_2=(0,a,0)$, $a_3=(L,0,0)$ and

\[
\widetilde\Gamma_{G_2}=\langle t_{a_1},\,t_{a_3},\,\sigma_{D,a_2/2}\rangle,\qquad G_2=\mathbb R^3/\widetilde\Gamma_{G_2},
\]

\noindent where $\sigma_{D,a_2/2}(x,y,z)=(-x,y+\frac a2,-z)$ and $\sigma_{D,a_2/2}^{\,2}=t_{a_2}$. The element
$g=t_{a_3}\circ\sigma_{D,a_2/2}\colon(x,y,z)\mapsto(L-x,\,y+\frac a2,\,-z)$ exchanges the boundary planes of $S$,
and it is the half-turn about the line $\{x=\frac L2,\,z=0\}$ composed with the translation by $\frac{a_2}{2}$. The group $\Gamma_S$ is generated by $t_{a_1}$ and $g$, whose restrictions to the central plane are

\[
\beta(y,z)=(y,z+b),\qquad \alpha(y,z)=(y+\tfrac a2,-z).
\]

Thus $\Gamma'_S=\langle\alpha,\beta\rangle$ is the fundamental group of the Klein bottle $K_{a,b}$, and
$\Gamma'_{S,0}=\langle\alpha^2,\beta\rangle$ is its translation lattice, which are exactly the data of $\TK_{a,b,d}$
with $d=L/4$. Hence

\[
G_2\cong\mathrm{Double}\bigl(\TK_{a,b,L/4}\bigr).
\]

\item \textbf{The double of $\KK_{a,b,d}$.} Let $P=\{x=0\}$, $D=\{x=z=0\}$, and

\[
a_1=(0,a,0),\qquad a_2=(0,0,2b),\qquad a_3=(L,0,0).
\]

Put

\[
\widetilde\Gamma_{B_3}=\langle\sigma_{D,a_1/2},\,\sigma_{P,a_2/2},\,t_{a_3}\rangle,\qquad B_3=\mathbb R^3/\widetilde\Gamma_{B_3},
\]

\noindent with $\sigma_{D,a_1/2}(x,y,z)=(-x,y+\frac a2,-z)$, $\sigma_{P,a_2/2}(x,y,z)=(-x,y,z+b)$,
$\sigma_{D,a_1/2}^{\,2}=t_{a_1}$ and $\sigma_{P,a_2/2}^{\,2}=t_{a_2}$. The holonomy group is
$\mathbb Z_2\times\mathbb Z_2$. The elements

\[
g_1=t_{a_3}\circ\sigma_{D,a_1/2}\colon(x,y,z)\mapsto(L-x,\,y+\tfrac a2,\,-z),\qquad
g_2=t_{a_3}\circ\sigma_{P,a_2/2}\colon(x,y,z)\mapsto(L-x,\,y,\,z+b)
\]

\noindent both exchange the boundary planes of $S$, and their restrictions to the central plane are

\[
\alpha(y,z)=(y+\tfrac a2,-z),\qquad \beta(y,z)=(y,z+b).
\]

Hence $\Gamma'_S=\langle\alpha,\beta\rangle=\pi_1(K_{a,b})$, and, since both generators exchange the two
sides, $\Gamma'_{S,0}=\langle\alpha\beta,\beta^2\rangle$. In Section~\ref{classificationdim3} the manifold
$\KK_{a,b,d}$ corresponds instead to $\Gamma'_0=\langle\alpha,\beta^2\rangle$. The two pairs
$(\Gamma',\Gamma'_0)$ are conjugate by the translation $c\colon z\mapsto z+\frac b2$: indeed
$c\alpha c^{-1}=\beta\alpha$, so $c\langle\alpha,\beta^2\rangle c^{-1}=\langle\beta\alpha,\beta^2\rangle
=\langle\alpha\beta,\beta^2\rangle$, while $c$ preserves $\langle\alpha,\beta\rangle$. Conjugating by the
translation $(x,y,z)\mapsto(x,y,z+\frac b2)$ therefore identifies $S/\Gamma_S$ with $\KK_{a,b,L/4}$, and

\[
B_3\cong\mathrm{Double}\bigl(\KK_{a,b,L/4}\bigr).
\]

\noindent Explicitly, with $s_1(x,y,z)=(-x,y,z+b)$ and $s_2(x,y,z)=(x,y+\frac a2,-z)$ as in
Section~\ref{classificationdim3}, the conjugation by $c$ gives $c\,\sigma_{P,a_2/2}\,c^{-1}=s_1$,
$c\,\sigma_{D,a_1/2}\sigma_{P,a_2/2}\,c^{-1}=s_2$ and $c\,t_{a_3}c^{-1}=t_{a_3}$, so
$c\,\widetilde\Gamma_{B_3}\,c^{-1}=\langle s_1,s_2,t_{a_3}\rangle$.

\end{itemize}

In each case $L=4d$, and the translation lattice of $\widetilde\Gamma$ contains $t_{(4d,0,0)}$; this is the
normalization used in the next subsections, where $d=d_i=L_i/4$ varies along the sequence.

\begin{lemma}\label{lem:proj}

Let $\hat\Gamma$ be a discrete group of isometries of $\mathbb R^m\times\mathbb R$ whose elements all have the form
$(g,\theta)$ with $g\in\mathrm{Isom}(\mathbb R^m)$, $\theta\in\mathrm{Isom}(\mathbb R)$, and such that $(\mathrm{id},t_\ell)\in\hat\Gamma$ for some $\ell>0$. Let $G$ be the image of $\hat\Gamma$ in $\mathrm{Isom}(\mathbb R^m)$ and assume $G$ is discrete. Then
$d_{GH}(\mathbb R^{m+1}/\hat\Gamma,\mathbb R^m/G)\le\ell/4$.

\end{lemma}

\begin{proof}

Write $p=(u,z)$ and $q=(u',z')$ for points of $\mathbb R^m\times\mathbb R$, and let $\bar d([u],[u'])=\min_{g\in G}|u-gu'|$ be the quotient metric of $\mathbb R^m/G$. The map $[(u,z)]\mapsto[u]$ is well defined and onto. For $(g,\theta)\in\hat\Gamma$, $|p-(g,\theta)q|\ge|u-gu'|$, so $d([p],[q])\ge\bar d([u],[u'])$. Conversely, take $g$ realizing $\bar d$ and a lift $(g,\theta)$. Composing with $(\mathrm{id},t_{k\ell})$ we may assume $|z-\theta z'|\le\ell/2$, so $d([p],[q])\le\bar d+\ell/2$. The graph has distortion $\le\ell/2$.

\end{proof}

\subsection{Proof of Theorem B}

Let $(M_i,g_i)$ be as in Theorem~B. By Section~\ref{classificationdim3}, each $M_i$ is isometric to a manifold of one of the five families
$T^2\times[-d,d]$, $K\times[-d,d]$, $\TT_{v,v',d}$, $\TK_{a,b,d}$, $\KK_{a,b,d}$. After passing to a subsequence, all the $M_i$ belong to the same family.
\medskip

Within each family, invertible linear maps of the covering slabs conjugate the displayed group actions as the parameters vary. They descend to diffeomorphisms onto a fixed compact manifold $M$; pulling back
the metrics by these diffeomorphisms, we may assume $M_i=(M,g_i)$, and Theorem~A applies. We treat the five families in turn. Throughout,
$\epsilon_i\to0$ denotes the parameter of Gromov--Hausdorff approximations.

\medskip

\noindent\emph{How the limits are computed.} In each of the five items below, Theorem~A is used only for the
dichotomy $\tau_\infty=0$ or $\tau_\infty>0$ and for the identification of the limit in case (ii). When
$\tau_\infty>0$, the limit $X$ is computed directly: we write down an explicit compact flat orbifold $O$, an
isometric involution $\bar\iota_0$ of $O$, and Gromov--Hausdorff approximations $f_i\colon\widetilde M_i\to O$ with
$f_i\circ\iota_i=\bar\iota_0\circ f_i$ up to an error tending to $0$, and we check that $M_i$ converges to
$O/\langle\bar\iota_0\rangle$ by an explicit approximation (the map $F_i$ in the two product cases, and
Lemma~\ref{lem:proj} applied to the groups $\hat\Gamma_i$ in the other three). Taking $\varphi_i(\iota_i)=\bar\iota_0$ and
$\psi_i(\bar\iota_0)=\iota_i$, the maps $f_i$ (together with quasi-inverses of the $f_i$) are equivariant approximations from
$(\widetilde M_i,\langle\iota_i\rangle)$ to $(O,\langle\bar\iota_0\rangle)$. These data agree with those of
Theorem~A. Indeed, let $(O_A,\bar\iota)$ and $f_i^A\colon\widetilde M_i\to O_A$ be the limit and the
approximations of Theorem~\ref{equivariant} and Lemma~\ref{lem:involution_limit} along the same subsequence,
and let $\bar f_i^A\colon O_A\to\widetilde M_i$ be quasi-inverses of the $f_i^A$. By
Lemma~\ref{lem:involution_limit}(i), $\bar f_i^A\circ\bar\iota$ is uniformly close to $\iota_i\circ\bar f_i^A$, so
the maps $\Phi_i=f_i\circ\bar f_i^A\colon O_A\to O$ are Gromov--Hausdorff approximations with error tending to $0$
and $\sup_{y\in O_A}d_O\bigl(\Phi_i(\bar\iota y),\bar\iota_0\Phi_i(y)\bigr)\to0$. A standard Arzel\`a--Ascoli
argument on a countable dense subset of $O_A$ gives a subsequence converging to an isometry
$\Phi\colon O_A\to O$ with $\Phi\circ\bar\iota=\bar\iota_0\circ\Phi$. We may therefore identify $(O_A,\bar\iota)$
with $(O,\bar\iota_0)$, and in the items below we write $\bar\iota$ for $\bar\iota_0$. In particular the
description of $X$ and the statement that $\bar\iota$ is a reflection of $O$ do not depend on the choice of
approximations in Theorem~\ref{equivariant}.

\medskip

\begin{itemize}

\item \textbf{The case $M\cong T^2\times[-d,d]$.} Here $M_i=(T_i\times[-d_i,d_i],g_i)$ with $T_i=\mathbb R^2/\Lambda_i$ a flat torus and $g_i$ the product metric. The boundary is $T_i\times\{\pm d_i\}$, disconnected, each component isometric to $T_i$, and $\tau_i=d_i$. Write $S^1_\ell=\mathbb R/\ell\mathbb Z$. Choose a Gauss-reduced basis $e_1^i,e_2^i$ of $\Lambda_i$, meaning $|e_1^i|\le|e_2^i|$ and $|\langle e_1^i,e_2^i\rangle|\le|e_1^i|^2/2$ (such a basis exists for every lattice in $\mathbb R^2$ by
Lagrange--Gauss reduction, see \cite{cassels}). Then $|\cos\vartheta_i|\le|e_1^i|/(2|e_2^i|)\le\frac12$, so the angle $\vartheta_i$ between them lies in $[\pi/3,2\pi/3]$, and we put $h_i=|e_2^i|\sin\vartheta_i\ge\frac{\sqrt3}{2}|e_2^i|$. Then $\operatorname{Area}(T_i)\ge\frac{\sqrt3}{2}|e_1^i||e_2^i|$, $\tfrac{\sqrt3}{4}|e_2^i|\le\tfrac{h_i}{2}\le\operatorname{diam}(T_i)\le|e_2^i|$, and $\operatorname{diam}(M_i)^2=\operatorname{diam}(T_i)^2+4d_i^2$. So $d_i\le D/2$ and $|e_2^i|\le 4D/\sqrt3$, and after passing to a subsequence $d_i\to d_\infty=\tau_\infty$, $|e_1^i|\to\lambda_1$, $|e_2^i|\to\lambda_2$, $h_i\to h_\infty$.

\bigskip

\noindent\emph{Limits of $T_i$.} If $\lambda_1>0$, after applying orthogonal changes of coordinates and taking a further subsequence, the reduced bases converge to a basis of a lattice $\Lambda_\infty$ and $T_i\to T_\infty=\mathbb R^2/\Lambda_\infty$. If $\lambda_1=0<\lambda_2$, let $u_i\perp e_1^i$ be a unit vector. Then $\langle\Lambda_i,u_i\rangle=h_i\mathbb Z$, and $[p]\mapsto\langle p,u_i\rangle \bmod h_i$ is a well defined onto map $T_i\to S^1_{h_i}$. Its fibres are circles of length $|e_1^i|$, and one checks that it changes distances by at most $|e_1^i|/2$. Composing with the homothety $S^1_{h_i}\to S^1_{h_\infty}$ (factor $\to1$) gives $\epsilon_i$-approximations $q_i:T_i\to S^1_{h_\infty}$, so $T_i\to S^1_{h_\infty}$. If $\lambda_2=0$, then $\operatorname{diam}(T_i)\to0$ and $T_i\to\{\mathrm{pt}\}$. Let $C_\infty$ denote this limit ($S^1_{h_\infty}$ or a point) when $\lambda_1=0$, and $q_i:T_i\to C_\infty$ the approximations.

\bigskip

\noindent\emph{Case $\tau_\infty=0$.} By Theorem~A(ii)(a), $X=\lim_{GH}T_i$, since both boundary components are isometric to $T_i$. Hence $X$ is a flat torus $T_\infty$ if $\lambda_1>0$, the circle $S^1_{h_\infty}$ if $\lambda_1=0<\lambda_2$, and a point if $\lambda_2=0$. Here $\bar\iota=\mathrm{Id}$ and $X=O$.

\bigskip

\noindent\emph{Case $\tau_\infty>0$.} Since $\operatorname{Vol}(M_i)=2d_i\operatorname{Area}(T_i)\to0$ and $d_i\to d_\infty>0$, we get $\operatorname{Area}(T_i)\to0$, so $|e_1^i|^2\le\frac{2}{\sqrt3}\operatorname{Area}(T_i)\to0$, i.e.\ $\lambda_1=0$. The double is $\widetilde M_i=T_i\times S^1_{4d_i}$, with $\iota_i(p,s)=(p,2d_i-s)$ (the two copies of $[-d_i,d_i]$ glued along their endpoints form the circle $S^1_{4d_i}$). Put $\phi_i(s)=\frac{d_\infty}{d_i}s$, $O=C_\infty\times S^1_{4d_\infty}$ and $f_i=q_i\times\phi_i$. Since $|\sqrt{a'^2+b'^2}-\sqrt{a^2+b^2}|\le|a'-a|+|b'-b|$ and $\phi_i$ has distortion $\le|\frac{d_\infty}{d_i}-1|\,2d_i\to0$, the map $f_i$ is onto with distortion tending to $0$. Hence $\widetilde M_i\to O$, a rectangular flat torus if $\lambda_2>0$ and a circle if $\lambda_2=0$, consistent with \cite{betderpic}, Theorem~D. With $\bar\iota(c,s)=(c,2d_\infty-s)$ we have $f_i\circ\iota_i=\bar\iota\circ f_i$ exactly, so $\bar\iota$ is the limit of $\iota_i$ in the sense of Lemma~\ref{lem:involution_limit}(i). It is non-trivial, in agreement with Lemma~\ref{lem:trivial_action}, and it is a pure reflection of the circle factor with $\operatorname{Fix}(\bar\iota)=C_\infty\times\{\pm d_\infty\}$. Therefore

\[
X=O/\langle\bar\iota\rangle=C_\infty\times[-d_\infty,d_\infty],
\]

\noindent that is, the annulus $S^1_{h_\infty}\times[-d_\infty,d_\infty]$ with two mirror circles (the orbifold $**$) if $\lambda_2>0$, and the segment $[-d_\infty,d_\infty]$ if $\lambda_2=0$. (This agrees with the direct map $F_i(p,x)=(q_i(p),\frac{d_\infty}{d_i}x)$, which is an approximation $F_i\colon M_i\to C_\infty\times[-d_\infty,d_\infty]$.)

\bigskip

\noindent\emph{Conclusion.} The collapsed limits of $T^2\times[-d,d]$ are exactly a point, $S^1$, $T^2$ (with $\bar\iota=\mathrm{Id}$), and $I$, $S^1\times I$ with mirror boundary (with $\bar\iota$ a reflection). Every prescribed flat $2$-torus occurs by keeping its lattice fixed and taking $d_i\to0$. For explicit examples of all the listed types, use rectangular tori $\mathbb R^2/(\delta_i\mathbb Z\oplus\lambda_i\mathbb Z)$: take $d_i\to0$ with $(\delta_i,\lambda_i)$ constant, $(\delta_i\to0,\lambda)$, or $(\delta_i,\lambda_i)\to(0,0)$ to get $T^2$, $S^1_\lambda$, a point; and take $d_i=d>0$ with $(\delta_i\to0,\lambda)$ or $(\delta_i,\lambda_i)\to(0,0)$ to get $S^1_\lambda\times[-d,d]$ and $[-d,d]$. In particular, the Klein bottle, the M\"obius band and the pillowcase $S^2(2,2,2,2)$ do not occur. When $\tau_\infty>0$, the limiting doubling involution is the specified reflection of the circle factor.
\medskip

\item \textbf{The case $M\cong K\times[-d,d]$.} Here $M_i=(K_i\times[-d_i,d_i],g_i)$ with $g_i$ the product metric and $K_i=K_{a_i,b_i}=\mathbb R^2/\langle\alpha_i,\beta_i\rangle$ a flat Klein bottle, where $\alpha_i(y,z)=(y+\frac{a_i}2,-z)$ and $\beta_i(y,z)=(y,z+b_i)$. The boundary is $K_i\times\{\pm d_i\}$, which is disconnected, each component is isometric to $K_i$, and $\tau_i=d_i$. Write $S^1_\ell=\mathbb R/\ell\mathbb Z$. We have $\operatorname{Area}(K_i)=\frac{a_ib_i}{2}$, so $\operatorname{Vol}(M_i)=d_ia_ib_i$, and $\operatorname{diam}(M_i)^2=\operatorname{diam}(K_i)^2+4d_i^2$. Moreover

\[
\max\bigl(\tfrac{a_i}{4},\tfrac{b_i}{2}\bigr)\le\operatorname{diam}(K_i)\le\sqrt{\tfrac{a_i^2}{16}+\tfrac{b_i^2}{4}}.
\]

\noindent The lower bound holds because the maps $K_i\to S^1_{a_i/2}$, $[(y,z)]\mapsto y$, and $K_i\to[0,\frac{b_i}2]=\mathbb R/\langle z\mapsto -z,\,z\mapsto z+b_i\rangle$, $[(y,z)]\mapsto z$, are $1$-Lipschitz and onto. For the upper bound, given two points one can use a power of $\alpha_i$ to make the $y$-distance at most $\frac{a_i}4$, and the $z$-distance is at most $\frac{b_i}2$ for either sign. Hence $d_i\le D/2$, $a_i\le4D$, $b_i\le2D$, and after passing to a subsequence $d_i\to d_\infty=\tau_\infty$, $a_i\to a_\infty$, $b_i\to b_\infty$.

\bigskip

\noindent\emph{Limits of $K_i$.} If $a_\infty,b_\infty>0$, the groups $\langle\alpha_i,\beta_i\rangle$ converge to $\langle\alpha_\infty,\beta_\infty\rangle$ and $K_i\to K_\infty=K_{a_\infty,b_\infty}$, a flat Klein bottle. If $a_\infty>0=b_\infty$, the map $[(y,z)]\mapsto y$ is an onto map $K_i\to S^1_{a_i/2}$ whose fibres have diameter at most $\frac{b_i}2$. It changes distances by at most $\frac{b_i}{2}$, so $K_i\to S^1_{a_\infty/2}$. If $a_\infty=0<b_\infty$, the map $[(y,z)]\mapsto z$ is an onto map $K_i\to[0,\frac{b_i}2]$ that changes distances by at most $\frac{a_i}{2}$, so $K_i\to[0,\frac{b_\infty}2]$. If $a_\infty=b_\infty=0$, then $\operatorname{diam}(K_i)\to0$ and $K_i\to\{\mathrm{pt}\}$. When $a_\infty b_\infty=0$, let $C_\infty$ denote this limit ($S^1_{a_\infty/2}$, $[0,\frac{b_\infty}2]$ or a point) and $q_i:K_i\to C_\infty$ the corresponding $\epsilon_i$-approximations (composed with the obvious homothety of factor $\to1$).

\bigskip

\noindent\emph{Case $\tau_\infty=0$.} By Theorem~A(ii)(a), $X=\lim_{GH}K_i$, since both boundary components are isometric to $K_i$. Hence $X$ is the Klein bottle $K_\infty$ if $a_\infty,b_\infty>0$, the circle $S^1_{a_\infty/2}$ if $a_\infty>0=b_\infty$, the interval $[0,\frac{b_\infty}2]$ if $a_\infty=0<b_\infty$, and a point if $a_\infty=b_\infty=0$. Here $\bar\iota=\mathrm{Id}$ and $X=O$.

\bigskip

\noindent\emph{Case $\tau_\infty>0$.} Since $\operatorname{Vol}(M_i)=d_ia_ib_i\to0$ and $d_i\to d_\infty>0$, we get $a_ib_i\to0$, i.e.\ $a_\infty b_\infty=0$, so $C_\infty$ and $q_i$ are defined. The double is $\widetilde M_i=K_i\times S^1_{4d_i}$ with $\iota_i(p,s)=(p,2d_i-s)$, since the two copies of $[-d_i,d_i]$ glued along their endpoints form the circle $S^1_{4d_i}$. Put $\phi_i(s)=\frac{d_\infty}{d_i}s$, $O=C_\infty\times S^1_{4d_\infty}$ and $f_i=q_i\times\phi_i$. Since $|\sqrt{a'^2+b'^2}-\sqrt{a^2+b^2}|\le|a'-a|+|b'-b|$ and $\phi_i$ has distortion at most $|\frac{d_\infty}{d_i}-1|\,2d_i\to0$, the map $f_i$ is onto with distortion tending to $0$. Hence $\widetilde M_i\to O$, which is a rectangular flat torus $S^1_{a_\infty/2}\times S^1_{4d_\infty}$ if $a_\infty>0=b_\infty$, the cylinder $[0,\frac{b_\infty}2]\times S^1_{4d_\infty}$ (the orbifold $**$) if $a_\infty=0<b_\infty$, and a circle if $a_\infty=b_\infty=0$. This is consistent with \cite{betderpic}. With $\bar\iota(c,s)=(c,2d_\infty-s)$ we have $f_i\circ\iota_i=\bar\iota\circ f_i$ exactly, so $\bar\iota$ is the limit of $\iota_i$ in the sense of Lemma~\ref{lem:involution_limit}(i). It is non-trivial, in agreement with Lemma~\ref{lem:trivial_action}, and it is a pure reflection of the circle factor with $\operatorname{Fix}(\bar\iota)=C_\infty\times\{\pm d_\infty\}$. Therefore

\[
X=O/\langle\bar\iota\rangle=C_\infty\times[-d_\infty,d_\infty],
\]

\noindent that is, the annulus $S^1_{a_\infty/2}\times[-d_\infty,d_\infty]$ with two mirror circles (the orbifold $**$) if $a_\infty>0=b_\infty$, the rectangle $[0,\frac{b_\infty}2]\times[-d_\infty,d_\infty]$ with all four sides mirrors and four corner reflectors of order $2$ (the orbifold $*2222$) if $a_\infty=0<b_\infty$, and the segment $[-d_\infty,d_\infty]$ if $a_\infty=b_\infty=0$. (This agrees with the direct map $F_i(p,x)=(q_i(p),\frac{d_\infty}{d_i}x)$, which is an approximation $F_i\colon M_i\to C_\infty\times[-d_\infty,d_\infty]$.)

\bigskip

\noindent\emph{Conclusion.} The collapsed limits of $K\times[-d,d]$ are exactly a point, $S^1$, $I$ and the Klein bottle $K_\infty$ (with $\bar\iota=\mathrm{Id}$), and $I$, $S^1\times I$ with mirror boundary and $I\times I$ with mirror boundary (with $\bar\iota$ a reflection). Each occurs. Take $d_i\to0$ with $(a_i,b_i)$ constant, $(a,b_i\to0)$, $(a_i\to0,b)$, or $(a_i,b_i)\to(0,0)$ to get $K_{a,b}$, $S^1_{a/2}$, $[0,\frac b2]$, a point. Take $d_i=d>0$ with $(a,b_i\to0)$, $(a_i\to0,b)$, or $(a_i,b_i)\to(0,0)$ to get $S^1_{a/2}\times[-d,d]$, $[0,\frac b2]\times[-d,d]$, $[-d,d]$. In particular a torus $T^2$ occurs as a limit of $K\times[-d,d]$ for no choice of parameters. The pillowcase and the M\"obius band do not occur. When $\tau_\infty>0$, the quotient by the limiting doubling involution is one of the displayed interval, annulus or rectangle limits; the Klein bottle occurs when $\tau_\infty=0$ and $\bar\iota=\mathrm{Id}$.
\medskip

\item \textbf{The case $M\cong \TT_{v,v',d}$.} Here $M_i=([-d_i,d_i]\times\mathbb R^2)/\Gamma_i$, where $\Gamma_i$ is generated by $t_{v_i}$ and $\varsigma_i(x,p)=(-x,p+\frac{v_i'}2)$. Let $\Lambda_i'=\mathbb Zv_i+\mathbb Z\frac{v_i'}2$ be the lattice of the soul $N_i=\mathbb R^2/\Lambda_i'$, and $\Lambda_{i,0}=\mathbb Zv_i+\mathbb Zv_i'$ its sublattice of index $2$. Let $\chi_i:\Lambda_i'\to\{\pm1\}$ be the character with kernel $\Lambda_{i,0}$. Then $\Gamma_i=\{(x,p)\mapsto(\chi_i(\mu)x,p+\mu):\mu\in\Lambda_i'\}$, the boundary is connected, $\tau_i=d_i$, and $\operatorname{Vol}(M_i)=2d_i\operatorname{Area}(N_i)$. Write $S^1_\ell=\mathbb R/\ell\mathbb Z$. Choose a Gauss-reduced basis $e_1^i,e_2^i$ of $\Lambda_i'$ as above (in particular $|e_1^i|\le|e_2^i|$ and $\vartheta_i\in[\pi/3,2\pi/3]$), and let $h_i=|e_2^i|\sin\vartheta_i\ge\frac{\sqrt3}2|e_2^i|$ be the height of the lattice over $e_1^i$, so that $\langle\Lambda_i',u_i\rangle=h_i\mathbb Z$ for a unit vector $u_i\perp e_1^i$. The projection $M_i\to N_i$ is $1$-Lipschitz, so $\operatorname{diam}(N_i)\le D$, and $d_i\le D$. Projection onto the height circle gives $\operatorname{diam}(N_i)\ge h_i/2\ge\frac{\sqrt3}{4}|e_2^i|$, so $|e_2^i|\le4D/\sqrt3$. After passing to a subsequence, $d_i\to d_\infty=\tau_\infty$, $|e_1^i|\to\lambda_1$, $|e_2^i|\to\lambda_2$, $h_i\to h_\infty$, and $\chi_i(e_1^i)$, $\chi_i(e_2^i)$ are constant. Since $\chi_i\ne1$, exactly one of the following holds:

\bigskip

\noindent (A) $\chi_i(e_1^i)=1,\ \chi_i(e_2^i)=-1$ (the short vector acts as a pure translation on the universal cover);  

\noindent (B) $\chi_i(e_1^i)=-1$ (the short vector acts as a glide, i.e.\ $(x,p)\mapsto(-x,p+e_1^i)$ is in $\Gamma_i$ and its square is $t_{2e_1^i}$).

\bigskip

\noindent The two cases give isometric manifolds up to a change of basis of $\Lambda_i'$ (the three index-$2$ sublattices of
$\Lambda_i'$ are permuted by $GL_2(\mathbb Z)$), but they do not give the same degeneration when $|e_1^i|\to0$,
because the collapsing direction $e_1^i$ is then distinguished. In case (A) the orientation reversal $x\mapsto-x$ is
carried only by the long direction $e_2^i$, so it survives in the limit as a glide along the circle that remains. In
case (B) it is carried by the short vector $e_1^i$, so as $e_1^i\to0$ it converges to the pure reflection
$x\mapsto-x$, and the twist disappears into a mirror. This is why the limits in case $\tau_\infty>0$ differ: a
M\"obius band in case (A) and an annulus in case (B).

\bigskip

\noindent\emph{Limits of $N_i$.} As for $T^2\times[-d,d]$: if $\lambda_1>0$, orthogonal changes of coordinates and a further subsequence give a nondegenerate limiting basis, and $N_i\to N_\infty=\mathbb R^2/\Lambda_\infty$ is a flat torus; if $\lambda_1=0<\lambda_2$ then $N_i\to S^1_{h_\infty}$ via $[p]\mapsto\langle p,u_i\rangle\bmod h_i$, which changes distances by at most $|e_1^i|/2$; if $\lambda_2=0$ then $N_i\to\{\mathrm{pt}\}$.

\bigskip

\noindent\emph{Case $\tau_\infty=0$.} The boundary is connected, and by Theorem~\ref{theorem:slab}(4b) $\partial M_i/\sigma_i\cong N_i$. By Theorem~A(ii)(b), $X=\lim N_i$, that is, a flat torus $N_\infty$ if $\lambda_1>0$, the circle $S^1_{h_\infty}$ if $\lambda_1=0<\lambda_2$, and a point if $\lambda_2=0$. Here $\bar\iota=\mathrm{Id}$ and $X=O$. (Note that $\lim\partial M_i$ itself may differ from $X$, since $\partial M_i=\mathbb R^2/\Lambda_{i,0}$ has the longer lattice.)

\bigskip

\noindent\emph{Case $\tau_\infty>0$.} Since $2d_i\operatorname{Area}(N_i)\to0$ and $d_\infty>0$, $\operatorname{Area}(N_i)\to0$, so $|e_1^i|^2\le\frac2{\sqrt3}\operatorname{Area}(N_i)\to0$, i.e.\ $\lambda_1=0$. The double is $\widetilde M_i=\mathbb R^3/\widetilde\Gamma_i$ with

\[
\widetilde\Gamma_i=\{(x,p)\mapsto(\chi_i(\mu)x+4d_ik,\ p+\mu):\ \mu\in\Lambda_i',\ k\in\mathbb Z\}
\]

\noindent (the group generated by $\Gamma_i$ and $t_{(4d_i,0)}$), and $\iota_i$ is induced by $\rho_i(x,p)=(2d_i-x,p)$. This group acts freely, since an element with $\chi_i(\mu)=-1$ has $\mu\ne0$. Write $\mathbb R^3=\mathbb R^2_{(x,u)}\times\mathbb R e_1^i$ (with $u$ the coordinate along $u_i$). Every element of $\widetilde\Gamma_i$, and of $\hat\Gamma_i:=\langle\widetilde\Gamma_i,\rho_i\rangle$ (so that $M_i=\mathbb R^3/\hat\Gamma_i$), is a product of an isometry of $\mathbb R^2_{(x,u)}$ and a translation of $\mathbb Re_1^i$. Let $Q_i$ and $\hat Q_i=\langle Q_i,\rho_i'\rangle$ be the images of $\widetilde\Gamma_i$ and $\hat\Gamma_i$ in $\mathrm{Isom}(\mathbb R^2_{(x,u)})$, where $\rho_i'(x,u)=(2d_i-x,u)$. Put $\kappa_i(x,u)=(x+4d_i,u)$. In case (A), $t_{e_1^i}\in\widetilde\Gamma_i$ acts trivially on $\mathbb R^2_{(x,u)}$ and translates $\mathbb Re_1^i$ by $\ell_i=|e_1^i|$, and $Q_i=\langle\gamma_i,\kappa_i\rangle$ with $\gamma_i(x,u)=(-x,u+h_i)$. In case (B), the elements of $\widetilde\Gamma_i$ acting trivially on $\mathbb R^2_{(x,u)}$ are the $t_{2me_1^i}$, so $\ell_i=2|e_1^i|$, and, since for each value $n$ of $\langle\mu,u_i\rangle/h_i$ both signs of $\chi_i$ occur (whatever $\chi_i(e_2^i)$ is, change the parity of the $e_1^i$-coefficient), $Q_i=\langle r,\kappa_i,(x,u)\mapsto(x,u+h_i)\rangle$ with $r(x,u)=(-x,u)$. Both $\widetilde\Gamma_i$ and $\hat\Gamma_i$ contain the translation $t_{\ell_ie_1^i/|e_1^i|}$ along $\mathbb Re_1^i$ acting trivially on $\mathbb R^2_{(x,u)}$, so by Lemma~\ref{lem:proj}

\[
d_{GH}(\widetilde M_i,\mathbb R^2/Q_i)\le\ell_i/4\to0,\qquad d_{GH}(M_i,\mathbb R^2/\hat Q_i)\le\ell_i/4\to0 .
\]

If $\lambda_2>0$, then $h_\infty>0$. The linear map $(x,u)\mapsto(d_\infty x/d_i,h_\infty u/h_i)$ conjugates $Q_i$ to the group $Q_\infty$ with the same generators and limiting parameters, and intertwines $\rho_i'$ with $\rho_\infty'$. Its two scale factors tend to $1$, so the induced maps of the compact quotients have distortion tending to $0$. Together with the projection estimate, this gives $\widetilde M_i\to O=\mathbb R^2/Q_\infty$. If $\lambda_2=0$, then $h_i\to0$ and one applies the lemma once more to the translation $\gamma_i^2=(0,2h_i)$ in case (A), or $(0,h_i)$ in case (B), to get $O=\mathbb R/\langle r,\kappa\rangle=[0,2d_\infty]$. The involution $\bar\iota$ is induced by $\rho_\infty'(x,u)=(2d_\infty-x,u)$ (all the approximations commute with $\rho_i$ up to the factor $d_\infty/d_i\to1$), so it is non-trivial, in agreement with Lemma~\ref{lem:trivial_action}, and it is a reflection with $\operatorname{Fix}(\bar\iota)$ the image of $\{x=\pm d_\infty\}$, the ambient limit of the boundary images in the sense of Remark~\ref{rem:boundary_limits}. Identifying $O$ and $X=O/\langle\bar\iota\rangle$ gives:

\begin{itemize}

\item[(A), $\lambda_2>0$:] $O=\mathbb R^2/\langle\gamma_\infty,\kappa_\infty\rangle$ is the flat Klein bottle $K_{2h_\infty,4d_\infty}$ (glide $\gamma_\infty$ with shift $h_\infty$, translation $\kappa_\infty$ by $4d_\infty$). The group $\langle\gamma_\infty,\kappa_\infty,\rho_\infty'\rangle$ is generated by the mirrors $x=\pm d_\infty$ and the glide $\gamma_\infty$, which exchanges them, so $X=([-d_\infty,d_\infty]\times\mathbb R)/\gamma_\infty$ is the flat M\"obius band of width $2d_\infty$ and core length $h_\infty$, whose boundary circle (of length $2h_\infty$) is a mirror (the orbifold $*\times$).

\item[(B), $\lambda_2>0$:] $O=[0,2d_\infty]\times S^1_{h_\infty}$, an annulus with two mirror circles ($**$), and $X=[0,d_\infty]\times S^1_{h_\infty}$, again of type $**$: the mirror $\{x=0\}$ comes from the soul (the collapsed glide) and $\{x=d_\infty\}$ from $\partial M$.

\item[$\lambda_2=0$:] $O=[0,2d_\infty]$ and $X=[0,d_\infty]$.

\end{itemize}

These agree with the direct picture: in case (A), $M_i$ tends to $([-d_\infty,d_\infty]\times\mathbb R)/\gamma_\infty$, and in case (B), to $([-d_\infty,d_\infty]/\pm)\times S^1_{h_\infty}$.

\bigskip

\noindent\emph{Conclusion.} The collapsed limits of $\TT_{v,v',d}$ are exactly a point, $S^1$, $T^2$ (with $\bar\iota=\mathrm{Id}$), and $[0,d]$, the annulus $S^1\times[0,d]$ with two mirror circles, and the M\"obius band with mirror boundary (with $\bar\iota$ a reflection). Every prescribed flat $2$-torus occurs by fixing a basis $(v,v'/2)$ of its lattice and taking $d_i\to0$. For explicit examples of all the listed types, take $v_i,v_i'$ orthogonal. For $d_i\to0$, the soul lattice $\mathbb Zv_i+\mathbb Z\frac{v_i'}2$ with fixed lengths, with one length $\to0$, or with both lengths $\to0$ gives $T^2$, $S^1$, a point. For $d_i=d>0$: $v_i=(\delta_i,0)$, $v_i'=(0,2h)$ with $\delta_i\to0$ gives the M\"obius band (case (A)) since the short vector $v_i$ has $\chi_i(v_i)=+1$; $v_i'=(2\delta_i,0)$, $v_i=(0,h)$ with $\delta_i\to0$ gives $[0,d]\times S^1_h$ (case (B)) since the short vector $v_i'/2$ has $\chi_i(v_i'/2)=-1$; and both lengths $\to0$ gives $[0,d]$. A M\"obius band occurs in this family.
\medskip

\item \textbf{The case $M\cong \KK_{a,b,d}$.} Here $M_i=([-d_i,d_i]\times\mathbb R^2)/\Gamma_i$ with $\Gamma_i=\langle s_1^i,s_2^i\rangle$,

\[
s_1^i(x,y,z)=(-x,y,z+b_i),\qquad s_2^i(x,y,z)=(x,y+\tfrac{a_i}2,-z).
\]

The soul is the Klein bottle $N_i=K_{a_i,b_i}=\mathbb R^2/\langle\alpha_i,\beta_i\rangle$, with $\alpha_i(y,z)=(y+\frac{a_i}2,-z)$ and $\beta_i(y,z)=(y,z+b_i)$. The boundary is the connected Klein bottle $K_{a_i,2b_i}$, $\tau_i=d_i$, $\operatorname{Area}(N_i)=\frac{a_ib_i}2$ and $\operatorname{Vol}(M_i)=d_ia_ib_i$. Write $S^1_\ell=\mathbb R/\ell\mathbb Z$. The maps $N_i\to S^1_{a_i/2}$, $[(y,z)]\mapsto y$, and $N_i\to[0,\frac{b_i}2]=\mathbb R/\langle z\mapsto-z,\,z\mapsto z+b_i\rangle$, $[(y,z)]\mapsto z$, are onto and $1$-Lipschitz, so $\operatorname{diam}(N_i)\ge\max(\frac{a_i}4,\frac{b_i}2)$. Since the projection $M_i\to N_i$ is $1$-Lipschitz, $\operatorname{diam}(N_i)\le D$ and $d_i\le D$. Hence $(d_i,a_i,b_i)$ is bounded, and after passing to a subsequence $d_i\to d_\infty=\tau_\infty$, $a_i\to a_\infty$, $b_i\to b_\infty$.

\bigskip

\noindent\emph{Limits of $N_i$.}

If $a_\infty,b_\infty>0$, the groups $\langle\alpha_i,\beta_i\rangle$ converge and $N_i\to K_\infty=K_{a_\infty,b_\infty}$.

If $a_\infty>0=b_\infty$, the map $[(y,z)]\mapsto y$ changes distances by at most $\frac{b_i}2$ (the $z$-difference can be reduced to at most $\frac{b_i}2$ using $\beta_i$), so $N_i\to S^1_{a_\infty/2}$.

If $a_\infty=0<b_\infty$, the map $[(y,z)]\mapsto z$ changes distances by at most $\frac{a_i}2$ (once the $z$-part is fixed, the $y$-difference can be reduced modulo $a_i$ using $\alpha_i^2$), so $N_i\to[0,\frac{b_\infty}2]$.

If $a_\infty=b_\infty=0$, then $\operatorname{diam}(N_i)\to0$ and $N_i\to\{\mathrm{pt}\}$.

\bigskip

\noindent\emph{Case $\tau_\infty=0$.} The boundary is connected and, by Theorem~\ref{theorem:slab}(4b), $\partial M_i/\sigma_i\cong N_i$. By Theorem~A(ii)(b), $X=\lim N_i$. Hence $X$ is the Klein bottle $K_\infty$ if $a_\infty,b_\infty>0$, the circle $S^1_{a_\infty/2}$ if $a_\infty>0=b_\infty$, the interval $[0,\frac{b_\infty}2]$ if $a_\infty=0<b_\infty$, and a point if $a_\infty=b_\infty=0$. Here $\bar\iota=\mathrm{Id}$ and $X=O$. (The limit of $\partial M_i=K_{a_i,2b_i}$ itself may differ from $X$.)

\bigskip

\noindent\emph{Case $\tau_\infty>0$.} Since $d_ia_ib_i\to0$ and $d_\infty>0$, we get $a_\infty b_\infty=0$. The double is $\widetilde M_i=\mathbb R^3/\widetilde\Gamma_i$ with $\widetilde\Gamma_i=\langle s_1^i,s_2^i,\theta_i\rangle$ and $\theta_i=t_{(4d_i,0,0)}$. This group is the conjugate of $\widetilde\Gamma_{B_3}$ by the translation $z\mapsto z+\frac{b_i}2$ (see Section~\ref{recollement}), so $\widetilde M_i$ is isometric to the manifold $B_3$ of that section. The involution $\iota_i$ is induced by $\rho_i(x,y,z)=(2d_i-x,y,z)$. Since $\rho_is_1^i=t_{(2d_i,0,b_i)}$, we have $M_i=\mathbb R^3/\hat\Gamma_i$ with $\hat\Gamma_i=\langle\widetilde\Gamma_i,\rho_i\rangle=\langle\rho_i,s_2^i,t_{(2d_i,0,b_i)}\rangle$. All elements of $\hat\Gamma_i$ act on each coordinate by an isometry of $\mathbb R$, so Lemma~\ref{lem:proj} applies along any coordinate axis, to $\widetilde\Gamma_i$ (giving $O$) and to $\hat\Gamma_i$ (giving $X$). It applies with $\ell=2b_i$ along $z$, since $(s_1^i)^2=t_{(0,0,2b_i)}$, and with $\ell=a_i$ along $y$, since $(s_2^i)^2=t_{(0,a_i,0)}$. When the remaining parameters converge to positive limits, the quotients of $\mathbb R^2$ by the image groups converge to the quotients by the limit groups, by a linear rescaling of the axes with factors tending to $1$. Write $\kappa_i(x)=x+4d_i$, $r_0(x,y)=(-x,y)$, and $\rho'_i(x,\cdot)=(2d_i-x,\cdot)$ for the induced maps.
\medskip

\emph{Sub-case $a_\infty>0=b_\infty$ (collapse of $z$).} The image of $\widetilde\Gamma_i$ in $\mathrm{Isom}(\mathbb R^2_{xy})$ is $\langle r_0,\ (x,y)\mapsto(x,y+\frac{a_i}2),\ \kappa_i\rangle$, so $d_{GH}(\widetilde M_i,[0,2d_i]\times S^1_{a_i/2})\le\frac{b_i}2\to0$ and $O=[0,2d_\infty]\times S^1_{a_\infty/2}$, an annulus with two mirror circles ($**$). Adding $\rho_i'$ gives $X=[0,d_\infty]\times S^1_{a_\infty/2}$, again of type $**$.
\medskip

\emph{Sub-case $a_\infty=0<b_\infty$ (collapse of $y$).} The image of $\widetilde\Gamma_i$ in $\mathrm{Isom}(\mathbb R^2_{xz})$ is $\widetilde G_i=\langle g_i,\varsigma,\kappa_i\rangle$, with $g_i(x,z)=(-x,z+b_i)$ and $\varsigma(x,z)=(x,-z)$, and $d_{GH}(\widetilde M_i,\mathbb R^2/\widetilde G_i)\le\frac{a_i}4\to0$. The group $\widetilde G_\infty$ contains the mirrors $z=kb_\infty$, the glide reflections $g_\infty$ along the perpendicular lines $x=2kd_\infty$, the rotations of order $2$ about the points $(2kd_\infty,(m+\frac12)b_\infty)$ (for instance $g_\infty\varsigma$ is the half-turn about $(0,\frac{b_\infty}2)$), and the rectangular lattice of translations generated by $(4d_\infty,0)$ and $(0,2b_\infty)$. So $O=\mathbb R^2/\widetilde G_\infty$ is of type $22{*}$ (pmg), of area $2d_\infty b_\infty$. Adding $\rho'_\infty$, i.e.\ the mirrors $x=(2k+1)d_\infty$, we get $\hat G_\infty=\langle\widetilde G_\infty,\rho'_\infty\rangle$. It has mirrors in both directions and the centred rectangular lattice generated by $(4d_\infty,0)$, $(0,2b_\infty)$ and $\rho'_\infty g_\infty=t_{(2d_\infty,b_\infty)}$, so $X=\mathbb R^2/\hat G_\infty$ is of type $2{*}22$ (cmm), of area $d_\infty b_\infty$. Concretely, $X$ is the rectangle $[0,d_\infty]\times[0,b_\infty]$ whose sides $x=d_\infty$, $z=0$, $z=b_\infty$ are mirrors, with two corner reflectors of order $2$ at $(d_\infty,0)$ and $(d_\infty,b_\infty)$, and whose side $x=0$ is folded at its midpoint $(0,\frac{b_\infty}2)$, a cone point of order $2$. The two corners $(0,0)$ and $(0,b_\infty)$ are identified and are not corners of the orbifold.
\medskip

\emph{Sub-case $a_\infty=b_\infty=0$.} Applying the lemma along $y$ and then along $z$ gives $O=[0,2d_\infty]$ and $X=[0,d_\infty]$.
\medskip

In every sub-case, $\bar\iota$ is induced by $\rho_\infty(x,\cdot)=(2d_\infty-x,\cdot)$, because the approximations in the lemma only forget the collapsed coordinates and commute with $\rho_i$ up to the factor $d_\infty/d_i\to1$. So $\bar\iota\neq\mathrm{Id}$, in agreement with Lemma~\ref{lem:trivial_action}, and it is a pure reflection in the $x$-coordinate whose fixed set is the image of $\{x=\pm d_\infty\}$, the ambient limit of the boundary images in the sense of Remark~\ref{rem:boundary_limits}. This also gives $X=O/\langle\bar\iota\rangle$, as in Theorem~A(i).

\bigskip

\noindent\emph{Conclusion.} The collapsed limits of $\KK_{a,b,d}$ are exactly the following.

\begin{center}

\renewcommand{\arraystretch}{1.2}

\begin{tabular}{|c|c|c|c|}

\hline
$d_\infty$ & $(a_\infty,b_\infty)$ & $O$ & $X$ \\ \hline
$0$ & $a_\infty,b_\infty>0$ & $K_{a_\infty,b_\infty}$ & $K_{a_\infty,b_\infty}$ \\
$0$ & $a_\infty>0=b_\infty$ & $S^1_{a_\infty/2}$ & $S^1_{a_\infty/2}$ \\
$0$ & $a_\infty=0<b_\infty$ & $[0,\frac{b_\infty}2]$ & $[0,\frac{b_\infty}2]$ \\
$0$ & $a_\infty=b_\infty=0$ & point & point \\ \hline
$>0$ & $a_\infty>0=b_\infty$ & $[0,2d_\infty]\times S^1_{a_\infty/2}$ $(**)$ & $[0,d_\infty]\times S^1_{a_\infty/2}$ $(**)$ \\
$>0$ & $a_\infty=0<b_\infty$ & pmg $(22{*})$ & cmm $(2{*}22)$ \\
$>0$ & $a_\infty=b_\infty=0$ & $[0,2d_\infty]$ & $[0,d_\infty]$ \\ \hline

\end{tabular}

\end{center}

\bigskip

Each occurs. Take $d_i\to0$ with $(a_i,b_i)=(a,b)$, $(a,b_i\to0)$, $(a_i\to0,b)$, $(a_i,b_i)\to(0,0)$, or $d_i=d>0$ with $(a,b_i\to0)$, $(a_i\to0,b)$, $(a_i,b_i)\to(0,0)$. In particular, the flat torus never occurs as a limit, and the only $2$-dimensional limits are the Klein bottle $(\times\times)$, $**$ and $2{*}22$.
\medskip

\item \textbf{The case $M\cong \TK_{a,b,d}$.} Here $M_i=([-d_i,d_i]\times\mathbb R^2)/\Gamma_i$ with $\Gamma_i=\langle s^i,t^i\rangle$,

\[
s^i(x,y,z)=(-x,\,y+\tfrac{a_i}2,\,-z),\qquad t^i(x,y,z)=(x,y,z+b_i).
\]

The soul is the Klein bottle $N_i=K_{a_i,b_i}=\mathbb R^2/\langle\alpha_i,\beta_i\rangle$, with $\alpha_i(y,z)=(y+\frac{a_i}2,-z)$ and
$\beta_i(y,z)=(y,z+b_i)$. The boundary is the connected rectangular torus $T_{a_i,b_i}=\mathbb R^2/\langle\alpha_i^2,\beta_i\rangle$, $\tau_i=d_i$, $\operatorname{Area}(N_i)=\frac{a_ib_i}2$ and
$\operatorname{Vol}(M_i)=d_ia_ib_i$. Write $S^1_\ell=\mathbb R/\ell\mathbb Z$. The maps $N_i\to S^1_{a_i/2}$, $[(y,z)]\mapsto y$, and $N_i\to[0,\frac{b_i}2]=\mathbb R/\langle z\mapsto-z,\,z\mapsto z+b_i\rangle$, $[(y,z)]\mapsto z$, are onto and $1$-Lipschitz, so $\operatorname{diam}(N_i)\ge\max(\frac{a_i}4,\frac{b_i}2)$. Since $M_i\to N_i$ is $1$-Lipschitz, $\operatorname{diam}(N_i)\le D$ and $d_i\le D$. Hence $(d_i,a_i,b_i)$ is bounded and, after passing to a subsequence, $d_i\to d_\infty=\tau_\infty$, $a_i\to a_\infty$, $b_i\to b_\infty$.

\bigskip

\noindent\emph{Limits of $N_i$.} The soul $N_i$ is the same Klein bottle $K_{a_i,b_i}$ as for $\KK_{a_i,b_i,d_i}$, so the limits are the same. If $a_\infty,b_\infty>0$, then $N_i\to K_\infty=K_{a_\infty,b_\infty}$. If $a_\infty>0=b_\infty$, the map $[(y,z)]\mapsto y$ changes distances by at most $\frac{b_i}2$, so $N_i\to S^1_{a_\infty/2}$. If $a_\infty=0<b_\infty$, the map $[(y,z)]\mapsto z$ changes distances by at most $\frac{a_i}{2}$ (reduce the $y$-difference modulo $a_i$ using $\alpha_i^2$), so $N_i\to[0,\frac{b_\infty}2]$. If $a_\infty=b_\infty=0$, then $N_i\to\{\mathrm{pt}\}$.

\bigskip

\noindent\emph{Case $\tau_\infty=0$.} The boundary is connected and, by Theorem~\ref{theorem:slab}(4b), $\partial M_i/\sigma_i\cong N_i$. By Theorem~A(ii)(b), $X=\lim N_i$. So $X$ is the Klein bottle $K_\infty$ if $a_\infty,b_\infty>0$, the circle $S^1_{a_\infty/2}$ if $a_\infty>0=b_\infty$, the interval $[0,\frac{b_\infty}2]$ if $a_\infty=0<b_\infty$, and a point if $a_\infty=b_\infty=0$. Here $\bar\iota=\mathrm{Id}$ and $X=O$. (The limit of $\partial M_i=T_{a_i,b_i}$ itself may differ from $X$; for instance, if $a_\infty,b_\infty>0$ it is a torus.)

\bigskip

\noindent\emph{Case $\tau_\infty>0$.} Since $d_ia_ib_i\to0$ and $d_\infty>0$, we get $a_\infty b_\infty=0$. The double is $\widetilde M_i=\mathbb R^3/\widetilde\Gamma_i$ with $\widetilde\Gamma_i=\langle s^i,t^i,\theta_i\rangle$, $\theta_i=t_{(4d_i,0,0)}$, that is the manifold $G_2$ of Section~\ref{recollement}, and $\iota_i$ is induced by $\rho_i(x,y,z)=(2d_i-x,y,z)$. One checks $\rho_is^i\rho_i^{-1}=\theta_i s^i$ and $\rho_it^i\rho_i^{-1}=t^i$. Put $w_i=\rho_is^i\colon(x,y,z)\mapsto(x+2d_i,y+\frac{a_i}2,-z)$. Then $s^i=\rho_iw_i$ and $\theta_i=w_i^2(s^i)^{-2}$, so

\[
M_i=\mathbb R^3/\hat\Gamma_i,\qquad\hat\Gamma_i=\langle\widetilde\Gamma_i,\rho_i\rangle=\langle\rho_i,\,t^i,\,w_i\rangle .
\]

Every element of $\widetilde\Gamma_i$ and $\hat\Gamma_i$ acts on each coordinate by an isometry of $\mathbb R$, so Lemma~\ref{lem:proj} applies along any coordinate axis, to $\widetilde\Gamma_i$ (giving $O$) and to $\hat\Gamma_i$ (giving $X$). Along $z$ it applies with $\ell=b_i$, since $t^i=t_{(0,0,b_i)}\in\Gamma_i$, and along $y$ with $\ell=a_i$, since $(s^i)^2=t_{(0,a_i,0)}$. When the remaining parameters converge to positive limits, the quotients of $\mathbb R^2$ by the image groups converge to the quotient by the limit group, by a linear rescaling of the axes with factors tending to $1$. Write $\kappa_i(x,\cdot)=(x+4d_i,\cdot)$ for the induced translation.
\medskip

\emph{Sub-case $a_\infty>0=b_\infty$ (collapse of $z$).} The image of $\widetilde\Gamma_i$ in $\mathrm{Isom}(\mathbb R^2_{xy})$ is $\widetilde G_i=\langle g_i,\kappa_i\rangle$ with $g_i(x,y)=(-x,y+\frac{a_i}2)$, because $t^i$ maps to the identity. Then $d_{GH}(\widetilde M_i,\mathbb R^2/\widetilde G_i)\le\frac{b_i}4\to0$, and $O=\mathbb R^2/\langle g_\infty,\kappa_\infty\rangle$ is a flat Klein bottle $K_{a_\infty,4d_\infty}$ (a glide reflection along $y$ of shift $\frac{a_\infty}2$, with a translation $4d_\infty$ perpendicular to its axis; it has no singular points). Adding $\rho'_\infty(x,y)=(2d_\infty-x,y)$, the group $\hat G_\infty=\langle\rho'_\infty,g_\infty,\kappa_\infty\rangle$ has mirrors $x=d_\infty+2kd_\infty$ and glide reflections $\rho'_\infty\circ(\text{translation by }(2d_\infty,\frac{a_\infty}2))=g_\infty$ along the parallel axes $x=2kd_\infty$, and a centred rectangular translation lattice generated by $(\pm2d_\infty,\frac{a_\infty}2)$. So $X=([-d_\infty,d_\infty]\times\mathbb R)/g_\infty$ is a M\"obius band of width $2d_\infty$ and core length $\frac{a_\infty}2$, whose boundary circle (of length $a_\infty$) is a mirror (the orbifold $*\times$). Its area is $d_\infty a_\infty=\frac12\operatorname{Area}(O)$.
\medskip

\emph{Sub-case $a_\infty=0<b_\infty$ (collapse of $y$).} The image of $\widetilde\Gamma_i$ in $\mathrm{Isom}(\mathbb R^2_{xz})$ is $\widetilde G_i=\langle r,\,t_{(0,b_i)},\,\kappa_i\rangle$, where $r(x,z)=(-x,-z)$ is the image of $s^i$, a half-turn, and $d_{GH}(\widetilde M_i,\mathbb R^2/\widetilde G_i)\le\frac{a_i}4\to0$. The group $\widetilde G_\infty$ is generated by the half-turn $r$ and the rectangular lattice $\langle(4d_\infty,0),(0,b_\infty)\rangle$, so $O=\mathbb R^2/\widetilde G_\infty$ is the pillowcase $S^2(2,2,2,2)$ (type $2222$, $p2$), of area $2d_\infty b_\infty$. Adding $\rho'_\infty(x,z)=(2d_\infty-x,z)$ we get $\hat G_\infty=\langle\widetilde G_\infty,\rho'_\infty\rangle$. Its translation lattice is still $\langle(4d_\infty,0),(0,b_\infty)\rangle$ and its linear parts form $\mathbb Z_2\times\mathbb Z_2$. It contains the mirrors $x=d_\infty+2kd_\infty$, the half-turns about the points $(2kd_\infty,\frac m2b_\infty)$, and the glide reflections $\rho'_\infty r\colon(x,z)\mapsto(x+2d_\infty,-z)$ along the horizontal axes $z=\frac m2b_\infty$ with shift $2d_\infty$. It has no reflection in a horizontal axis, since every element with linear part $\mathrm{diag}(1,-1)$ has $x$-shift in $2d_\infty+4d_\infty\mathbb Z$. So $X=\mathbb R^2/\hat G_\infty$ is of type $22{*}$ ($pmg$), of area $d_\infty b_\infty=\frac12\operatorname{Area}(O)$. Concretely, $X=([-d_\infty,d_\infty]\times S^1_{b_\infty})/\langle(x,z)\mapsto(-x,-z)\rangle$ is a disc with two cone points of order $2$, at $(0,0)$ and $(0,\frac{b_\infty}2)$, and with mirror boundary (the image of $x=\pm d_\infty$).
\medskip

\emph{Sub-case $a_\infty=b_\infty=0$.} Applying the lemma along $y$ and then along $z$ gives $O=[0,2d_\infty]$ and $X=[0,d_\infty]$.
\medskip

In every sub-case $\bar\iota$ is induced by $\rho_\infty(x,\cdot)=(2d_\infty-x,\cdot)$, since the approximations of the lemma only forget the collapsed coordinates and commute with $\rho_i$ up to the factor $d_\infty/d_i\to1$. Hence $\bar\iota\ne\mathrm{Id}$, in agreement with Lemma~\ref{lem:trivial_action}, and it is a pure reflection in the $x$-coordinate whose fixed set is the image of $\{x=\pm d_\infty\}$, the ambient limit of the boundary images in the sense of Remark~\ref{rem:boundary_limits}. So $X=O/\langle\bar\iota\rangle$, as in Theorem~A(i).

\bigskip

\noindent\emph{Conclusion.} The collapsed limits of $\TK_{a,b,d}$ are exactly the following.

\begin{center}

\renewcommand{\arraystretch}{1.2}

\begin{tabular}{|c|c|c|c|}

\hline
$d_\infty$ & $(a_\infty,b_\infty)$ & $O$ & $X$ \\ \hline
$0$ & $a_\infty,b_\infty>0$ & $K_{a_\infty,b_\infty}$ & $K_{a_\infty,b_\infty}$ \\
$0$ & $a_\infty>0=b_\infty$ & $S^1_{a_\infty/2}$ & $S^1_{a_\infty/2}$ \\
$0$ & $a_\infty=0<b_\infty$ & $[0,\frac{b_\infty}2]$ & $[0,\frac{b_\infty}2]$ \\
$0$ & $a_\infty=b_\infty=0$ & point & point \\ \hline
$>0$ & $a_\infty>0=b_\infty$ & $K_{a_\infty,4d_\infty}$ $(\times\times)$ & M\"obius band, mirror boundary $(*\times)$ \\
$>0$ & $a_\infty=0<b_\infty$ & $S^2(2,2,2,2)$ $(2222)$ & $22{*}$ \\
$>0$ & $a_\infty=b_\infty=0$ & $[0,2d_\infty]$ & $[0,d_\infty]$ \\ \hline

\end{tabular}

\end{center}

\bigskip

Each occurs: take $d_i\to0$ with $(a_i,b_i)=(a,b)$, $(a,b_i\to0)$, $(a_i\to0,b)$, $(a_i,b_i)\to(0,0)$, or $d_i=d>0$ with $(a,b_i\to0)$, $(a_i\to0,b)$, $(a_i,b_i)\to(0,0)$. In particular the flat torus never occurs as a limit, and the $2$-dimensional limits are $K_\infty$ $(\times\times)$, $*\times$ and $22{*}$.

\end{itemize}

\bigskip

\noindent\textbf{End of the proof of Theorem~B.} By the five items, the possible limits are the following.

\begin{center}

\renewcommand{\arraystretch}{1.2}

\begin{tabular}{|c|c|c|}

\hline
Family & $\tau_\infty=0$ & $\tau_\infty>0$ \\ \hline
$T^2\times[-d,d]$ & point, $S^1$, $T^2$ & $I$, $**$ \\
$K\times[-d,d]$ & point, $S^1$, $I$, $K$ & $I$, $**$, $*2222$ \\
$\TT_{v,v',d}$ & point, $S^1$, $T^2$ & $I$, $**$, $*\times$ \\
$\KK_{a,b,d}$ & point, $S^1$, $I$, $K$ & $I$, $**$, $2{*}22$ \\
$\TK_{a,b,d}$ & point, $S^1$, $I$, $K$ & $I$, $*\times$, $22{*}$ \\ \hline

\end{tabular}

\end{center}

\bigskip

The union of the columns $\tau_\infty>0$ is $\{I,\ **,\ *2222,\ *\times,\ 2{*}22,\ 22{*}\}$, and the union of the columns $\tau_\infty=0$ is
$\{\text{point},\ S^1,\ I,\ T^2,\ K\}$; this is exactly the list of Theorem~B. In the case $\tau_\infty>0$, $\bar\iota$ is a non-trivial reflection
whose fixed set is the ambient limit of the boundary images as in Remark~\ref{rem:boundary_limits}, and every limit in the table is realized by the explicit sequences given in the ``Conclusion'' of each item.

\qed

\medskip

\noindent\textbf{Acknowledgments.} The author thanks the Lebanese Association for Scientific Research for hospitality and support.

\medskip

\noindent\textbf{AI disclosure.} Large language models were used during the preparation of this article to check arguments for gaps, verify formulas and inequalities, suggest and draft some arguments, and assist with the exposition. The main ideas of the proofs are the author's, and all arguments were checked by the author, who takes full responsibility for the content and correctness of the article.


\begin{thebibliography}{99}

\bibitem{betderpic}

Bettiol, R.\,G.; Derdzinski, A.; Piccione, P.; Teichm\"uller theory and collapse of flat manifolds. Annali di Matematica Pura ed Applicata 197 (2018), 1247--1268.


\bibitem{bieb1}

Bieberbach, L.; \"Uber die Bewegungsgruppen der Euklidischen R\"aume, I. Math. Ann. 70 (1911), 297--336.


\bibitem{bieb2}

Bieberbach, L.; \"Uber die Bewegungsgruppen der Euklidischen R\"aume, II. Die Gruppen mit einem endlichen Fundamentalbereich. Math. Ann. 72 (1912), 400--412.


\bibitem{boilmailport} 

  Boileau, M.; Maillot, S.; Porti, J.; \emph{Three-Dimensional Orbifolds and Their Geometric Structures}, Panoramas et Synth\`eses 15, Soci\'et\'e Math\'ematique de France, Paris, 2003. 

\bibitem{bbi}

Burago, D.; Burago, Y.; Ivanov, S.; A Course in Metric Geometry. Graduate Studies in Mathematics 33, American Mathematical Society, Providence RI, 2001.


\bibitem{bgp}
Burago, Y.\,D.; Gromov, M.\,L.; Perel'man, G.\,Y.; A.\,D. Alexandrov spaces with curvature bounded below. Russian Math. Surveys 47 (1992), no. 2, 1--58.


\bibitem{buser1}

Buser, P.; \emph{Geometry and Spectra of Compact {R}iemann Surfaces}. Progress in Mathematics, Birkh{\"a}user, Boston, 1992.


\bibitem{cassels}

Cassels, J.\,W.\,S.; \emph{An Introduction to the Geometry of Numbers}. Grundlehren der Mathematischen Wissenschaften 99, Springer-Verlag, Berlin, 1959.


\bibitem{charlap}

Charlap, L.\,S.; \emph{Bieberbach Groups and Flat Manifolds}. Universitext, Springer, New York, 1986.


\bibitem{cff}
Cheeger, J.; Fukaya, K.; Gromov, M.; Nilpotent structures and invariant metrics on collapsed manifolds. J. Amer. Math. Soc. 5 (1992), 327--372.


\bibitem{cheegergromov}
Cheeger, J.; Gromov, M.; Collapsing Riemannian manifolds while keeping their curvature bounded. I. J. Differential Geom. 23 (1986), 309--346.


\bibitem{cheegergromov2}
Cheeger, J.; Gromov, M.; Collapsing Riemannian manifolds while keeping their curvature bounded. II. J. Differential Geom. 32 (1990), 269--298.


\bibitem{conway}

Conway, J.\,H.; The orbifold notation for surface groups. In: \emph{Groups, Combinatorics and Geometry} (Durham, 1990), LMS Lecture Note Series 165, Cambridge Univ. Press, Cambridge, 1992, pp. 438--447.

\bibitem{davis}

Davis, M.\,W.; Lectures on orbifolds and reflection groups. In: \emph{Transformation Groups and Moduli Spaces of Curves}, ALM 16, Int. Press, 2011, pp. 63--93.

\bibitem{elaf}
El Mir, C.; Lafontaine, J.; Higher dimensional M\"obius bands and their boundaries. Manuscripta Math. 160 (2019), 187--198.


\bibitem{fp} Frigerio, R.; Petronio, C.; Construction and recognition of hyperbolic $3$-manifolds with geodesic boundary. Trans. Amer. Math. Soc. 356 (2004), 3243--3282.


\bibitem{fukaya1}

Fukaya, K.; Theory of convergence for Riemannian orbifolds. Jpn. J. Math. (N.S.) 12 (1986), 121--160.


\bibitem{fukaya}
Fukaya, K.; Collapsing Riemannian manifolds to ones of lower dimensions. J. Differential Geom. 25 (1987), 139--156.


\bibitem{fukayayama}

Fukaya, K.; Yamaguchi, T.; The fundamental groups of almost non-negatively curved manifolds. Ann. of Math. (2) 136 (1992), 253--333.


\bibitem{gromovalmost}
Gromov, M.; Almost flat manifolds. J. Differential Geom. 13 (1978), 231--241.


\bibitem{gromov}

Gromov, M.; \emph{Metric Structures for Riemannian and Non-Riemannian Spaces}. Modern Birkh\"auser Classics, Birkh\"auser Boston, 2007 (English edition).


\bibitem{hantzschewendt}

Hantzsche, W.; Wendt, H.; Dreidimensionale euklidische Raumformen. Math. Ann. 110 (1935), 593--611.

\bibitem{hy}

Huang, Z.; Yamaguchi, T.; Inradius collapsed manifolds with a lower Ricci curvature bound. Preprint, arXiv:2511.22081 (2025).

\bibitem{kang3d}

Kang, E.\,S.; Moduli spaces of 3-dimensional flat manifolds. J. Korean Math. Soc. 43 (2006), 1065--1080.


\bibitem{kangkim}

Kang, E.\,S.; Kim, J.\,Y.; Deformation spaces of 3-dimensional flat manifolds. Commun. Korean Math. Soc. 18 (2003), 95--104.


\bibitem{lange}
Lange, C.; Orbifolds from a metric viewpoint. Geom. Dedicata 209 (2020), 43--57.


\bibitem{nakajima}

Nakajima, S.; \"Uber konvexe Kurven und Fl\"achen. Tohoku Math. J. 29 (1928), 227--230.

\bibitem{pansu85}

Pansu, P.; Effondrement des vari\'et\'es riemanniennes, d'apr\`es J. Cheeger et M. Gromov. S\'eminaire Bourbaki 1983/84, expos\'e no. 618, Ast\'erisque 121--122 (1985), 63--82.


\bibitem{plesken}

Plesken, W.; Schulz, T.; Counting crystallographic groups in low dimensions. Exp. Math. 9 (2000), 407--411.


\bibitem{ratcliffe}

Ratcliffe, J.\,G.; \emph{Foundations of Hyperbolic Manifolds}, 2nd ed., Graduate Texts in Mathematics 149, Springer, 2006.

\bibitem{ruh}
Ruh, E.\,A.; Almost flat manifolds. J. Differential Geom. 17 (1982), 1--14.

\bibitem{scott}

Scott, P.; The geometries of $3$-manifolds. Bull. London Math. Soc. 15 (1983), no.~5, 401--487.


\bibitem{szczepanski1}

Szczepa\'nski, A.; Problems on Bieberbach groups and flat manifolds. Geom. Dedicata 120 (2006), 111--118.


\bibitem{szczepanski2}

Szczepa\'nski, A.; \emph{Geometry of Crystallographic Groups}, Vol. 4 of Algebra and Discrete Mathematics, World Scientific, 2012.


\bibitem{thurston2}
Thurston, W.\,P.; \emph{The Geometry and Topology of Three-Manifolds}, Princeton University lecture notes, 1978--1981.


\bibitem{tietze}

Tietze, H.; {\"U}ber Konvexheit im kleinen und im gro{\ss}en und {\"u}ber gewisse den Punkten einer Menge zugeordnete Dimensionszahlen. Math. Z. 28 (1928), 697--707.


\bibitem{wolfspaces}

Wolf, J.\,A.; \emph{Spaces of Constant Curvature}, 5th ed., Publish or Perish, 1984.

\bibitem{wong}

Wong, J.; An extension procedure for manifolds with boundary. Pacific J. Math. 235 (2008), 173--199.

\bibitem{yz}

Yamaguchi, T.; Zhang, Z.; Inradius collapsed manifolds. Geom. Topol. 23 (2019), 2793--2860.

\bibitem{yz1}

Yamaguchi, T.; Zhang, Z.; Limits of manifolds with boundary I. Preprint, arXiv:2406.00970 (2024).

\bibitem{yz2}

Yamaguchi, T.; Zhang, Z.; Limits of manifolds with boundary II. Preprint, arXiv:2504.05497 (2025).


\end{thebibliography}
\end{document}